\documentstyle[10pt,amscd,amssymb]{amsart}
\theoremstyle{plain}
 \newtheorem{thm}{Theorem}[section]
 \newtheorem{lem}[thm]{Lemma}
 \newtheorem{prop}[thm]{Proposition}
 \newtheorem{cor}[thm]{Corollary}
 
\theoremstyle{definition}
 \newtheorem{defn}{Definition}[section]
\theoremstyle{remark}
 \newtheorem{rem}{Remark}[section]

\newcommand{ \Supp}{\operatorname{Supp}}
\newcommand{\Ext}{\operatorname{Ext}}
\newcommand{\Hom}{\operatorname{Hom}}
\newcommand{\codim}{\operatorname{codim}}
\newcommand{\im}{\operatorname{im}}

\newcommand{\rk}{\operatorname{rk}}

\newcommand{\NS}{\operatorname{NS}}
\newcommand{\coker}{\operatorname{coker}}
\newcommand{\Pic}{\operatorname{Pic}}
\newcommand{\ch}{\operatorname{ch}}
\newcommand{\td}{\operatorname{td}}

\newcommand{\Hilb}{\operatorname{Hilb}}
\newcommand{\Quot}{\operatorname{Quot}}

\newcommand{\Amp}{\operatorname{Amp}}

\font\b=cmr10 scaled \magstep5
\def\bigzerou{\smash{\lower1.7ex\hbox{\b 0}}}

\numberwithin{equation}{section}
\begin{document}

\title
{Irreducibility of moduli spaces of vector bundles on K3 surfaces}
\author[K\={o}ta Yoshioka]{K\={o}ta Yoshioka\\
}
 \address{
Department of mathematics, Faculty of Science, Kobe University,
Kobe, 657, Japan}
\email{yoshioka@@math.kobe-u.ac.jp}

\begin{abstract}
In this paper,
we show the moduli spaces of stable sheaves on K3 surfaces 
are irreducible symplectic manifolds, if
the associated Mukai vectors are primitive.
More precisely, we show that
they are related to the Hilbert scheme of points. 
We also compute the period of these spaces. 
As an application of our result, we discuss Montonen-Olive
duality in Physics. In particular our computations of Euler
characteristics of moduli spaces are compatible with Physical
computations by Minahan et al.
\end{abstract}

 \maketitle

\section{Introduction}

\subsection{Main result}
Let $X$ be a projective K3 surface defined over ${\Bbb C}$
and $H$ an ample divisor on $X$.
Let $\omega$ be the fundamental class of $X$.
Let $E$ be a coherent sheaf on $X$.
By the identification $H^4(X,{\Bbb Z}) \cong {\Bbb Z} \omega$,
we regard the second Chern class $c_2(E)$ as an integer.
Since $(c_1(E)^2)$ is even, the second Chern character $\ch_2(E)$
belongs to ${\Bbb Z}$. 
We define the Mukai vector of $E$ by
\begin{equation}
\begin{split}
v(E):=&\ch(E)\sqrt{\td_X}\\
=&\rk(E)+c_1(E)+(\rk(E)+\ch_2(E))\omega \in H^*(X,{\Bbb Z}),
\end{split}
\end{equation}
where we identify $H^0(X,{\Bbb Z})$ with ${\Bbb Z}$
and $\td_X=1+2 \omega$ is the Todd class of $X$.
For an element $v \in H^*(X,{\Bbb Z})$,
we denote the $0$-th component $v_0 \in H^0(X,{\Bbb Z})$
by $\rk v$ and the second component $v_1 \in H^2(X,{\Bbb Z})$
by $c_1(v)$.
We set $\ell(v):=\gcd(\rk v,c_1(v)) \in {\Bbb Z}_{\geq 0}$.
Then $v$ is written as 
$v=\ell(v)(r+\xi)+a \omega$,
where $r \in {\Bbb Z}$, $\xi \in H^2(X,{\Bbb Z})$ and $r+\xi$ is primitive. 
We denote the moduli space of stable sheaves $E$ of
$v(E)=v$ by $M_H(v)$.
If $v$ is primitive and $H$ is general in the ample cone
$\Amp(X)$ of $X$
(i.e. there are at most countable number of 
hyperplanes $W_n \subset H^2(X,{\Bbb Q})$, 
$n=1,2,\dots$ which depends on $v$ and $H$ belongs to
$\Amp(X) \setminus \cup_n W_n$ [Y3]), then 
$M_H(v)$ is a smooth projective scheme.
In [Mu1], Mukai showed that $M_H(v)$ has a symplectic structure.
In order to get more precise information,
Mukai [Mu2] introduced a quite useful notion called
Mukai lattice $(H^*(X,{\Bbb Z}),\langle\;\;,\;\;\rangle)$,
where the pairing is defined by
\begin{equation}
\begin{split}
\langle x,y \rangle:=&-\int_X x^{\vee}y\\
=&\int_X(x_1y_1-x_0y_2-x_2y_0),
\end{split}
\end{equation}
where $x_i \in H^{2i}(X,{\Bbb Z})$ 
(resp. $y_i \in H^{2i}(X,{\Bbb Z})$) 
is the $2i$-th component of $x$ (resp. $y$)
and $x^{\vee}=x_0-x_1+x_2$.
Hence $\langle\;\;,\;\;\rangle$ is an integral primitive bilinear form on
$H^*(X,{\Bbb Z})$.
By the language of this lattice, 
we can write down Riemann-Roch theorem
in a simple form:
We set 
\begin{equation}
\chi(E,F):=\sum_{i=0}^2 \dim \Ext^i(E,F)
\end{equation}
for coherent sheaves $E$ and $F$.
Then Riemann-Roch theorem implies that
\begin{equation}
 \chi(E,F)=-\langle v(E),v(F) \rangle.
\end{equation}
In particular we get that $\dim M_H(v)=\langle v^2 \rangle+2$.
 
If $v$ is a primitive isotropic vector,
then $M_H(v)$ is a surface with a symplectic structure.
Mukai proved that $M_H(v)$ is a K3 surface and described the period
in terms of Mukai lattice. 
If $v$ is a primitive Mukai vector of $\langle v^2 \rangle >0$,
then $M_H(v)$ is a higher dimensional symplectic manifold.
If $\rk v=1$, then $M_H(v)$ is the Hilbert scheme of points on $X$.
Indeed every torsion free sheaf of rank 1 is give by
$I_Z \otimes L$, where $I_Z$ is the ideal sheaf of a 0-dimensional
subscheme of $X$ and $L$ is a line bundle of $c_1(L)=c_1(v)$. 
Beauville [B] proved that it is an example of higher dimensional 
irreducible symplectic manifold.
For an irreducible symplectic manifold,
Beauville [B] defined the period and proved local Torelli theorem.
As an example, he also computed the period of Hilbert scheme of points
on $X$.  
For higher rank cases,
Mukai [Mu3] (rank 2 case), O'Grady [O1] ($\ell(v)=1$ case)
and the author [Y5] ($\langle v^2 \rangle>2 \ell(v)^2$ or $\ell(v)=1$ case)
proved that $M_H(v)$ is an irreducible symplectic manifold
and described the period of $M_H(v)$ in terms of
Mukai lattice.
For classification of $M_H(v)$,
it is important to determine the period.
Indeed, it is a birational invariant ([Mu3]), and
affirmative solution of 
Torelli conjecture will imply that an irreducible symplectic manifold
is determined by its period, up to birational equivalence.

In this paper, by using [Y5] extensively,
we prove the following theorem,
which is expected by many people (for example, see [D], [Mu3], [O1]).
\begin{thm}\label{thm:main}
Let $v$ 
be a primitive Mukai vector such that $\rk v>0$ and 
$c_1(v) \in \NS(X)$.
\begin{enumerate}
\item[(1)]
$M_H(v)$ is not empty for a general ample divisor $H$
if and only if $\langle v^2 \rangle \geq -2$.
\item[(2)]
Assume that $\langle v^2 \rangle \geq -2$.
Then for a general ample divisor $H$,
\begin{enumerate}
\item[(2-1)]
$M_H(v)$ is obtained by compositions of deformations and birational
transformations from $\Hilb_X^{\langle v^2 \rangle/2+1}$.
In particular $M_H(v)$ is an irreducible symplectic manifold.
\item[(2-2)]
Let $B_{M_H(v)}$ be Beauville's bilinear form on $H^2(M_H(v),{\Bbb Z})$.
Then
$$
\theta_v:(v^{\perp},\langle\;\;\;,\;\;\; \rangle) \to 
(H^2(M_H(v),{\Bbb Z}),B_{M_H(v)})
$$
is an isometry which preserves Hodge structures
for $\langle v^2 \rangle \geq 2$,
where $\theta_v:v^{\perp} \to H^2(M_H(v),{\Bbb Z})$ 
is the canonical homomorphism
defined by using a quasi-universal family.

\end{enumerate}
\end{enumerate}
\end{thm}
Here we only use deformations of $M_H(v)$ induced by
deformation of complex structures of $X$.

Since birationally equivalent Calabi-Yau manifolds have the same 
Hodge numbers ([Ba],[De-L]),
we get the following Corollary.
\begin{cor}\label{cor:hodge}
Keep the notations as above.
Then
$h^{p,q}(M_H(v))=h^{p,q}(\Hilb_X^{\langle v^2 \rangle/2+1})$.
In particular, $\chi(M_H(v))=\chi(\Hilb_X^{\langle v^2 \rangle/2+1})$.
\end{cor}

In [V-W], Vafa and Witten considered a partition function 
$Z_r^{\alpha}(\tau)$, $\alpha \in H^2(X,{\Bbb Z})$,
$\tau \in {\Bbb H}:=\{z \in {\Bbb C}| \Im z>0 \}$
associated with $N=4$ super symmetric Yang-Mills theory on
a 4 manifold $X$.
Under suitable vanishing conditions (e.g. $H^0(X,ad(E) \otimes K_X)=0$),
it is related to ``Euler characteristics'' of moduli spaces
of vector bundles.
For a K3 surface case,
$Z_r^{\alpha}(\tau)$ is given by
\begin{equation}
 Z_r^{\alpha}(\tau)=
\sum
\begin{Sb}
\rk v=r\\
c_1(v)=\alpha
\end{Sb}
``\chi(M_H(v))"q^{\langle v^2 \rangle/2r},
\end{equation}
where $``\chi(M_H(v))"$ is a kind of ``Euler characteristics''
of a suitable compactification of $M_H(v)$.
Recently, this invariant was computed in [MNVW].
In section 4, by using Corollary \ref{cor:hodge},
we shall check that
their computation coincides with the Euler characteristics
of $M_H(v)$, if $v$ is primitive.

For a non-primitive Mukai vector, we have the following existence
condition.
\begin{cor}
Let $v$ be a Mukai vector of $\rk v>0$.
Then there is a semi-stable sheaf $E$ of $v(E)=v$
with respect to a general ample divisor $H$
if and only if $v=nw$, $n \in {\Bbb Z}$,
$w \in H^*(X,{\Bbb Z})$ with
$\langle w^2 \rangle \geq -2$.
\end{cor}

\subsection{Outline of the proof} 
We shall explain how to prove (2-1) of this theorem.
In [Y3], we discussed chamber structure of polarizations.
Let $v=l(r+\xi)+a \omega$, $\xi \in \NS(X)$
be a Mukai vector of $l=\ell(v)$ and $r>0$.
We choose an ample divisor $H$ on $X$
which does not lie on walls with respect to $v$.
Then
\begin{itemize}
\item[$(\natural)$]
for every $\mu$-semi-stable sheaf $E$ of $v(E)=v$,
if $F \subset E$ satisfies
$(c_1(F),H)/\rk F=(c_1(E),H)/\rk E$, then
$c_1(F)/\rk F=c_1(E)/\rk E$.
\end{itemize}
Thus $v(F)=l'(r+\xi)+a' \omega$ for some $l',a'$.
In particular, if $v$ is primitive, then 
$M_H(v)$ is compact.
Let ${\cal M}(v)$ be the stack of coherent sheaves $E$ of
$v(E)=v$.
We shall fix a general ample divisor $H$ with respect to $v$.
${\cal M}(v)^{\mu ss}$ (resp. ${\cal M}(v)^{\mu s}$) 
denotes the open substack of 
${\cal M}(v)$ consisting of
$\mu$-semi-stable sheaves (resp. $\mu$-stable sheaves).

In [Y5], we proved Theorem \ref{thm:main}
under the assumption $\langle v^2 \rangle>2l^2$ or $l=1$.
Hence we may assume that $\langle v^2 \rangle \leq 2l^2$ and $l >1$. 
However for convenience sake of the reader, 
we only use the results for $l=1$ case.
We note that isometry group ${\mathrm{O}}(H^*(X,{\Bbb Z}))$ of Mukai lattice
acts transitively on the set
$V_n:=\{x \in H^*(X,{\Bbb Z})|\text{$x$ is primitive}, 
\langle x^2 \rangle=2n \}$ and 
${\mathrm{O}}(H^*(X,{\Bbb Z}))/\pm 1$ is generated by the following 3 kinds of 
isometries:
\begin{enumerate}
\item
Translation:
For $N \in \Pic(X)$, 
\begin{equation}
\begin{matrix}
T_N: & H^*(X,{\Bbb Z}) & \to & H^*(X,{\Bbb Z})\\
& x & \mapsto & \ch(N) x
\end{matrix}
\end{equation}
is an isometry.
\item
${\mathrm{O}}(H^2(X,{\Bbb Z}))$ acts on ${\mathrm{O}}(H^*(X,{\Bbb Z}))$.
\item
Reflection:
For a $(-2)$- vector $v_1 \in H^*(X,{\Bbb Z})$,
\begin{equation}
\begin{matrix}
R_{v_1}: & H^*(X,{\Bbb Z}) & \to & H^*(X,{\Bbb Z})\\
& x & \mapsto & x+\langle x,v_1 \rangle v_1
\end{matrix}
\end{equation}
is an isometry.
\end{enumerate}
Therefore it is very important to understand reflections. 
\newline
Geometric realization of reflections:
As we shall see in Corollary \ref{cor:FMT}, a reflection is realized as a
Fourier-Mukai transform.
Here we shall explain a special case.   
Let $E_1$ be a stable vector bundle of $\Ext^1(E_1,E_1)=0$ ($E_1$ is called
exceptional vector bundle).
Since $\Ext^2(E_1,E_1) \cong \Hom(E_1,E_1)^{\vee} ={\Bbb C}$,
Riemann-Roch theorem implies that 
$\langle v(E_1),v(E_1) \rangle=-\chi(E_1,E_1)=-2$.
Thus $v_1:=v(E_1)$ is a $(-2)$-vector.
Let $E$ be a stable vector bundle of $v(E)=v$.
Assume that 
\begin{enumerate}
\item[(a)] 
$\Ext^i(E_1,E)=0$, $i=1,2$.
\item[(b)]
$\phi:E_1 \otimes \Hom(E_1,E) \to E$
is surjective and $\ker \phi$ is stable.
\end{enumerate}
Then $w:=v(\ker \phi)$ is given by $\chi(E_1,E)v_1-v=
-(v+\langle v,v_1 \rangle v_1)$.
Thus $-v(\ker \phi)$ is the $(-2)$-reflection of $v$ by $v_1$.
Hence under conditions (a) and (b), $(-2)$-reflection of
Mukai lattice induces a birational map
$M_H(v) \cdots \to M_H(w)$.
Replacing $\ker \phi$ by
$\coker(\phi^{\vee}:E^{\vee} \to E_1^{\vee}\otimes \Hom(E_1,E)^{\vee})$,
we may replace (b) by the condition (b'):
\begin{enumerate}
\item[(b')]
$\phi:E_1 \otimes \Hom(E_1,E) \to E$
is surjective in codimension 1 and $\ker \phi$ is stable.
\end{enumerate}
 Thus under (a) and (b'), we have a birational map
 $M_H(v) \cdots \to M_H(w^{\vee})$.
We would like to apply this story for a suitable pair of $v_1$ and $v$
which satisfy $\langle v_1, v \rangle=-1$.
In order to get condition (b'), 
we shall prove the following key lemma which was proved 
under the assumption
$lr<r_1<(l+1)r$ in [Y5, Prop. 4.5].
\begin{lem}\label{lem:-D*R}
Let $(X,H)$ be a polarized smooth projective surface of
$\NS(X)={\Bbb Z}H$.
Let $(r_1,d_1)$ and $(r,d)$ be pairs of integers
such that $r_1, r>0$ and $d r_1-r d_1=1$.
We assume that $lr<r_1$.
Let $E_1$ be a $\mu$-stable vector bundle of $\rk(E_1)=r_1$ and
$\deg(E_1)=d_1$, where $\deg(E_1)=(c_1(E_1),H)/(H^2)$.
\begin{enumerate}
\item[(1)]
Let $E$ be a $\mu$-stable sheaf of $\rk(E)=lr$ and 
$\deg(E)=ld$.
Then every non-zero homomorphism $\varphi:E_1 \to E$ is 
surjective in codimension 1 and $\ker \varphi$ is a $\mu$-stable sheaf.
\item[(2)]
Let $E'$ be a $\mu$-stable vector bundle of 
$\rk(E')=r_1-lr$ and $\deg(E')=d_1-ld$.
Let $\phi:E' \to E_1$ be a non-zero homomorphism.
Then $\phi$ is injective and $E:=\coker \phi$
is a $\mu$-semi-stable sheaf.
\end{enumerate}
\end{lem}
For the proof of this lemma,
we use the following fact:
\begin{itemize}
\item
We consider the triangle in ${\Bbb R}^2$ with vertices
$(0,0)$, $(r_1-lr,d_1-ld)$ and $(r_1,d_1)$.
Then there is no integral point in the interior of this 
triangle.
\end{itemize}
Indeed, this condition gives a strong restriction on homomorphisms
$\varphi$, $\phi$ and the Harder-Narasimhan
polygons of $\ker \varphi$ and $\coker \phi$.
The proof will be done in Preliminaries.

In order to use this lemma,
we need to compare ${\cal M}(v)^{\mu ss}$ and ${\cal M}(v)^{\mu s}$.
More precisely,
we need dimension counting of various constructible substacks of
${\cal M}(v)^{\mu ss}$.
Technically this is the most important part in this paper.
In [D-L], Drezet and Le Potier computed the dimension
of the substack of non-semi-stable sheaves.
In their computation,
the existence of exceptional vector bundle 
is very important.
In our case, we concentrate our consideration on 
${\cal M}(v)^{\mu ss}$.
By our assumption on $H$,
exceptional vector bundle $E$ of
$v(E)=r +\xi+b \omega$, $b \in {\Bbb Z}$
is important.
Hence we divide our proof into two cases:
\begin{enumerate}
\item[A.]
There is no $(-2)$ vector of the form
$r+\xi+b \omega$, i.e.
$((\xi^2)+2)/2r \not \in {\Bbb Z}$.
\item[B.]
There is a $(-2)$ vector of the form 
$r+\xi+b \omega$, i.e.
$((\xi^2)+2)/2r  \in {\Bbb Z}$.
\end{enumerate}

In section 2, we treat case A.
In particular, we prove the following inequality:
\begin{equation}\label{ineq:1}
\dim ({\cal M}(v)^{\mu ss} \setminus {\cal M}(v)^{\mu s})
\leq \langle v^2 \rangle.
\end{equation}
Then we can apply Lemma \ref{lem:-D*R}.
For suitable choice of 
\begin{enumerate}
\item
a primitive Mukai vector $v:=l(r+d H)+a \omega$ on $(X,H)$ and
\item
an exceptional vector bundle $E_1$ of $v(E_1)=r_1+d_1 H+a_1 \omega$,
\end{enumerate}
we can construct a birational map
$M_H(v) \cdots \to M_H(w^{\vee})$
sending a general $\mu$-stable vector bundle $E \in M_H(v)$
to $F:=\coker(E^{\vee} \to \Hom(E_1,E)^{\vee} \otimes E_1^{\vee})
 \in M_H(w^{\vee})$
where $H$ and $v'$ satisfy that
(1) $(H^2)/2=(r_1 \langle v^2 \rangle/2l+r)r_1/l-r^2>0$,
(2) $\langle v_1,v \rangle=-1$ and
(3) $\ell(w^{\vee})=1$ and hence Theorem \ref{thm:main} holds for 
$M_H(w^{\vee})$.
We remark that we need to choose {\it a sufficiently
large $r_1$} for the condition (1). 
In the same way as in [Y5, 4.3], we get Theorem \ref{thm:main}
for case A.
More precisely, considering deformations of $M_H(v)$ induced by 
deformations of $(X,H)$ and translations $T_N$,
we can reduce the problem to this situation.

In section 3,
we treat case B.
If $\langle v^2 \rangle \geq 2l^2$, then we also have the inequality
\eqref{ineq:1}, and hence the same proof as in case A works. 
If $\langle v^2 \rangle<2l^2$, then it is known that
there is no $\mu$-stable sheaf.
Hence we can not apply Lemma \ref{lem:-D*R} in this form.
When the $(-2)$ vector is $v({\cal O}_X)=1+\omega$,
T. Nakashima found the following fact:

We set $v=l-a \omega$. Then the inequality
$0 \leq \langle v^2 \rangle <2l^2$ implies that
$0 \leq a<l$.
We assume that $a \geq 2$.
Let $E$ be a $\mu$-stable vector bundle of
$v(E)=a-l \omega$.
Then $H^0(X,E)=H^2(X,E)=0$ and $\dim H^1(X,E)=l-a$.
We consider the universal extension (another example of reflection)
\begin{equation}
 0 \to E \to E' \to H^1(X,E) \otimes {\cal O}_X \to 0.
\end{equation}
It is easy to see that $E'$ is a stable vector bundle
and we get an immersion 
$M_H(a-l \omega)^{\mu s, loc} \hookrightarrow M_H(v)$,
where $M_H(a-l \omega)^{\mu s, loc}$ is the open subscheme of
$M_H(a-l \omega)$ consisting of $\mu$-stable vector bundles.

This result can be easily extended to general cases.
Hence what we should do is to prove the irreducibility
of $M_H(v)$ and the classification of $M_H(v)$
consisting of non-locally free sheaves.
By similar dimension counting as in case A, 
we shall classify non-locally free components
and prove the irreducibility of 
$M_H(v)$.
The classification of non-locally free components
of $M_H(v)$ is described as follows: 
\begin{prop}\label{prop:nlf}
Keep the notations in Theorem \ref{thm:main}.
Then 
$M_H(v)$ consists of non-locally free sheaves if and only if
$\rk v=1$, $v=(\rk v_0)v_0-\omega$ or $v=l -\omega$,
where $v_0$ is a Mukai vector of $\langle v_0^2 \rangle=-2$.
For these spaces, $M_H(v) \cong \Hilb_X^{\langle v^2 \rangle/2+1}$.
\end{prop}

\section{Preliminaries}

\subsection{Notation}
Except sections 1.5.1 and 1.5.2, we assume that $X$ is a K3 surface.
For a scheme $S$, $p_S:S \times X \to S$ denotes the projection.
For a Mukai vector $v$, we fix a general ample divisor $H$ which satisfies
$(\natural)$ in section 0.2.
Obviously for any subsheaf $E' \subset E$ of $\mu$-semi-stable sheaf
$E$ of $v(E)=v$,
if $c_1(E')/\rk E'=c_1(E)/\rk E$, then $v(E')$ also satisfies $(\natural)$.
${\cal M}(v), {\cal M}(v)^{\mu ss}$ and ${\cal M}(v)^{\mu s}$ are stacks in 
section 0.2.
${\cal M}(v)^{ss}$ and 
${\cal M}(v)^{s}$ 
denote the open substack of 
${\cal M}(v)$ consisting of
semi-stable sheaves and stable sheaves 
respectively.
Since ${\cal M}(v)^{\mu ss}$ is bounded,
it is a quotient stack of an open subscheme of 
some quot-scheme by some general linear group
(see Appendix).
Hence our dimension counting of substack of ${\cal M}(v)^{\mu ss}$
can be regarded as that of subscheme of some quot scheme.  
$M_H(v)^{\mu s}$ (resp. $M_H(v)^{loc}$) be 
the open subscheme of $M_H(v)$ consisting of
$\mu$-stable sheaves (resp. stable vector bundles).

Mukai homomorphism:
Let ${\cal E}$ be a quasi-universal family of similitude $\rho$ on
$M_H(v) \times X$, that is,
${\cal E}_{|\{E \} \times X} \cong E^{\oplus \rho}$
for all $E \in M_H(v)$ ([Mu3]).
By using ${\cal E}$, Mukai constructed a natural homomorphism
$$
\theta_v:v^{\perp} \longrightarrow H^2(M_H(v),{\Bbb Z})_f
$$
defined by
$$
\theta_v(x):=\frac{1}{\rho}
\left[p_{M_H(v)*}((\ch{\cal E})\sqrt{\td_X}x^{\vee})\right]_1,
$$
where $H^2(M_H(v),{\Bbb Z})_f$ is the torsion free quotient of
$H^2(M_H(v),{\Bbb Z})$. 
We note that $\theta_v$ does not 
depend on the choice of a quasi-universal family.

\subsection{Some results from [Y5]}
We collect some results which are necessary to prove 
Theorem \ref{thm:main}. 

\begin{thm}\label{thm:period}
Let $v=l(r+\xi)+a \omega,\; \xi \in H^2(X,{\Bbb Z})$ be a primitive
Mukai vector such that $l=\ell(v)$ and $r>0$.
If $\langle v^2 \rangle/2 \geq l^2$ or $l=1$, and
$H$ is general, then
$M_H(v)$ is obtained by compositions of deformations and birational 
transformations from $ \Hilb_X^{\langle v^2 \rangle/2+1}$.
In particular, $M_H(v)$ is an irreducible symplectic manifold.
Let $B_{M_H(v)}$ be Beauville's bilinear form on $H^2(M_H(v),{\Bbb Z})$.
If $\langle v^2 \rangle/2 > l^2$, or $l=1$ and $\langle v^2 \rangle/2=1$, 
then
$$
\theta_v:(v^{\perp},\langle\;\;\;,\;\;\; \rangle) \to 
(H^2(M_H(v),{\Bbb Z}),B_{M_H(v)})
$$
is an isometry which preserves Hodge structures
for $\langle v^2 \rangle \geq 2$.
\end{thm}
When $l=\ell(v)=1$, this theorem was first proved by O'Grady [O1].
In this paper, we only use this theorem for the case where $\ell(v)=1$.

The following is essentially due to O'Grady [O1].
We can see a different proof based on G\"{o}ttsche and Huybrechts' argument
[G-H] in [Y7]. 
\begin{prop}[{[Y5, Prop. 1.1]}]\label{prop:deform}
Let $X_1$ and $X_2$ be K3 surfaces, and let
$v_1:=l(r+\xi_1)+a_1 \omega \in H^*(X_1,{\Bbb Z})$
and $v_2:=l(r+\xi_2)+a_2 \omega \in H^*(X_2,{\Bbb Z})$
be primitive Mukai vectors such that
$(1)$ $r,l>0$, $(2)$ $r+\xi_1$ and $r+\xi_2$ are primitive,
$(3)$ $\langle v_1^2 \rangle=\langle v_2^2 \rangle=2s$,
and $(4)$ $a_1 \equiv a_2 \mod l$.
Then $M_{H_1}(v_1)$ and $M_{H_2}(v_2)$ are deformation equivalent.
In particular, $M_{H_1}(v_1)$ is an irreducible symplectic manifold
and $\theta_{v_1}$ is an isometry of Hodge structures if and only if 
$M_{H_2}(v_2)$ and $\theta_{v_2}$ have the same properties. 
\end{prop}

\begin{lem}[{[Y5, Lem. 5.1]}]\label{lem:inequality}
Let $x_1,x_2,x_3, y_1,y_2,y_3$ be integers such that $x_1,x_2,x_3>0$
and $y_1 x_3-x_1 y_3=1$.  
If 
\begin{equation}
\frac{y_1}{x_1}>\frac{y_2}{x_2}>\frac{y_3}{x_3},
\end{equation}
then $x_2 \geq x_1+x_3$.
\end{lem}

\begin{lem}[{[Y5, Lem. 4.1]}]\label{lem:4-1}
Let $(X,H)$ be a polarized smooth projective surface of
$\NS(X)={\Bbb Z}H$.
Let $(r_1,d_1)$ and $(r,d)$ be pairs of integers
such that $r_1, r>0$ and $d r_1-r d_1=1$.
Let $E_1$ be a $\mu$-stable vector bundle of $\rk(E_1)=r_1$ and
$\deg(E_1)=d_1$, where $\deg(E_1)=(c_1(E_1),H)/(H^2)$.
Let $E$ be a $\mu$-stable sheaf of $\rk(E)=lr$ and 
$\deg(E)=ld$.
Then the non-trivial extension
\begin{equation}
0 \to E_1 \to E' \to E \to 0
\end{equation}
is a $\mu$-stable sheaf.
\end{lem}

\begin{lem}[{[Y5, Lem. 4.4]}]\label{lem:quot}
Let $v$ be an arbitrary Mukai vector of $\rk v>0$.
Let ${\cal M}_H(v)^{\mu ss}$ be the moduli stack of 
$\mu$-semi-stable sheaves $E$ of $v(E)=v$, and
${\cal M}_H(v)^{p \mu ss}$ the closed substack of ${\cal M}_H(v)^{\mu ss}$ consisting
of properly $\mu$-semi-stable sheaves.
We assume that $\langle v^2 \rangle/2 \geq l^2$.
Then
\begin{equation}
\codim {\cal M}_H(v)^{p \mu ss} \geq \langle v^2 \rangle/2l-l+1.
\end{equation}
In particular, if ${\cal M}_H(v)^{\mu ss}$ is not empty,
then there is a $\mu$-stable sheaf $E$ of $v(E)=v$.   
\end{lem}

\begin{pf}
Since we need the proof of this lemma,
we shall give an outline of the proof.
For more details, see [Y5, sect. 5.3].
By Mukai [Mu1], we get that
\begin{equation}
\dim {\cal M}_H(v)^{\mu ss} \geq (\langle v^2 \rangle+2)-1. 
\end{equation}
We shall show that 
\begin{equation}
\dim {\cal M}_H(v)^{pss} \leq (\langle v^2 \rangle +1)-
(\langle v^2 \rangle/2l-l+1).
\end{equation}
For this purpose, we shall estimate the moduli number of
Jordan-H\"{o}lder filtrations. 
Let $E$ be a $\mu$-semi-stable sheaf of
$v(E)=v$ and let
$0 \subset F_1 \subset F_2 \subset \cdots \subset F_t=E$
be a Jordan-H\"{o}lder filtration of $E$ with respect to $\mu$-stability.
We set $E_i:=F_i/F_{i-1}$.
By using Lemma \ref{lem:A-est1} in Appendix successively,
we see that the moduli number of this filtration is bounded by
\begin{equation}
\sum_{i \leq j} (\dim \Ext^1(E_j,E_i)-\dim \Hom(E_j,E_i))
= -\chi(E,E)+\sum_{i>j}\chi(E_j,E_i)+
\sum_{i<j}\dim \Ext^2(E_j,E_i)+t.
\end{equation}
We set
$v(E):=lr+l\xi+a \omega$ and
$v(E_i):=l_i r+l_i \xi+a_i \omega$, where $\xi \in \NS(X)$.
Since $\langle v(E_i),v(E_j) \rangle=l_i l_j (\xi^2)-r(l_i a_j+l_j a_i)$,
we see that
\begin{equation}
\sum_{i>j}\chi(E_j,E_i) = -\sum_{i>j} \langle v(E_j),v(E_i) \rangle
=-\sum_{i}\frac{(l-l_i)\langle v(E_i)^2 \rangle}{2l_i}.
\end{equation}
We set $\max_i \{l_i\}=(l-k)$.
Let $i_0$ be an integer such that $ \langle v(E_{i_0})^2 \rangle \geq 0$.
Since $\sum_i l_i=l$, we obtain that $t \leq k+1$.  
Since $l-l_i-k \geq 0$ and $\langle v(E_i)^2 \rangle \geq -2$,
 we get that
\begin{align*}
\sum_{i>j} \langle v(E_j),v(E_i) \rangle 
&=k\sum_i \frac{\langle v(E_i)^2 \rangle}{2l_i}+
\sum_i \frac{(l-l_i-k)\langle v(E_i)^2 \rangle}{2l_i}\\
& \geq  k\frac{\langle v(E)^2 \rangle}{2l}-
\sum_{i \ne i_0} (l-l_i-k)\\
& \geq k\frac{\langle v(E)^2 \rangle}{2l}-(l-1-k)k.
\end{align*}
If $r>1$ or $l_i>1$ for some $i$,
then for a general filtration, there are $E_i$ and $E_j$ such that
$\Ext^2(E_j,E_i)=0$.
Therefore we get that 
$\sum_{i<j} \dim \Ext^2(E_j,E_i) \leq (k+1)k/2-1$
for a general filtration.
Then the moduli number of these filtrations is bounded by
\begin{align*}
\langle v^2 \rangle- k\frac{\langle v^2 \rangle}{2l}
+(l-1-k)k+\frac{(k+1)k}{2}-1+t 
&\leq 
(\langle v^2 \rangle+1)- \left(k\frac{\langle v^2 \rangle}{2l}
-lk+\frac{k(k-1)}{2}+1 \right)\\
& \leq (\langle v^2 \rangle+1)- \left(\frac{\langle v^2 \rangle}{2l}
-l+1 \right).
\end{align*}
Therefore we get a desired estimate for this case.
In particular, our lemma holds for the case where $\rk v/\ell(v)>1$.
For the case where $l_i=1$ for all $i$ and $r=1$, 
see [Y5, sect. 5.3].
\end{pf}

\subsection{Estimate on properly semi-stable sheaves}

\begin{lem}\label{lem:hom}
\begin{enumerate}
\item[(1)]
Let $E$ be a stable sheaf and $F$ a semi-stable sheaf
such that $v(F)/\rk F=v(E)/\rk E$.
Then $\Hom(E,F) \otimes E \to F$ is injective.
In particular $\dim \Hom(E,F) \leq \rk F/\rk E$.
\item[(2)] 
Let $E$ be a $\mu$-stable sheaf and $F$ a $\mu$-semi-stable sheaf
such that $c_1(F)/\rk F=c_1(E)/\rk E$.
Then $\Hom(E,F) \otimes E \to F$ is injective.
In particular $\dim \Hom(E,F) \leq \rk F/\rk E$.
\end{enumerate}
\end{lem}

\begin{lem}\label{lem:pss}
Let $v$ be a Mukai vector of $\langle v^2 \rangle>0$
(we don't assume the primitivity of $v$).
We set
\begin{equation}
{\cal M}(v)^{pss}:=\{E \in {\cal M}(v)^{ss}|
\text{$E$ is properly semi-stable }\}.
\end{equation}
Then $\dim {\cal M}(v)^{pss} \leq \langle v^2 \rangle$.
In particular,
if ${\cal M}(v)^{ss} \ne \emptyset$, then 
${\cal M}(v)^{s} \ne \emptyset$ and 
$\dim {\cal M}(v)^{ss}=\langle v^2 \rangle+1$.
\end{lem}
\begin{pf}
We set $v=lv'$, where $v'$ is a primitive Mukai vector.
We shall prove this lemma by induction on $l$.
Let $E_1$ be a stable sheaf of $v(E_1)=l_1v'$
and $E_2$ a semi-stable sheaf
of $v(E_2)=l_2v'$, where $l_1+l_2=l$.
By induction hypothesis,
$\dim {\cal M}^{ss}(v_i)=\langle v_i^2 \rangle+1$, $i=1,2$.
We shall estimate the dimension of
the substack $J(v_1,v_2)$ whose element $E$ fits in an extension
\begin{equation}
0 \to E_1 \to E \to E_2 \to 0.
\end{equation}
By Lemma \ref{lem:hom},
$\dim \Ext^2(E_2,E_1)=\dim \Hom(E_1,E_2) \leq l_2/l_1$.
Moreover if $E_1$ is general, then
$\Hom(E_1,E_2)=0$. Hence by Lemma \ref{lem:A-est1} in Appendix, we get 
\begin{equation}
\begin{split}
\dim J(v_1,v_2) & \leq
\dim {\cal M}(v_1)^{ss}+\dim {\cal M}(v_2)^{ss}+
 \langle v_1,v_2 \rangle+\max\{l_2/l_1-1,0 \}\\
 & =
\langle v_1^2 \rangle+\langle v_2^2 \rangle+2+
\langle v_1,v_2 \rangle+\max\{l_2/l_1-1,0\}\\
& =(\langle v^2 \rangle +1)-(\langle v_1,v_2 \rangle-\max\{l_2/l_1,1\}).
\end{split}
\end{equation}
Since 
\begin{equation}
\langle v_1,v_2 \rangle=l_1 \frac{\langle v_2^2 \rangle}{2l_2}+
l_2\frac{\langle v_1^2 \rangle}{2l_1} \geq l,
\end{equation}
we get $\dim J(v_1,v_2) \leq \langle v^2 \rangle$.
Therefore we get our lemma.
\end{pf}

\subsection{Semi-stable sheaves of isotropic Mukai vector}

\begin{lem}\label{lem:ss}
Let $w$ be a primitive isotropic Mukai vector of $\rk w>0$.
Then 
$\dim {\cal M}(lw)^{ss}=l$.
\end{lem}
\begin{pf}
Let $E$ be a semi-stable sheaf of $v(E)=lw$.
We shall first prove that
there are stable sheaves $E_1,E_2,\dots,E_k$ of $v(E_i)=w$
such that 
\begin{equation}\label{eq:JHF}
 E \cong \oplus_{i=1}^k F_i,
\end{equation}
where $F_i$ are $S$-equivalent to $E_i^{\oplus n_i}$.
\newline
Proof of the claim:
By the proof of Mukai [Mu2, Prop. 4.4]
(Fourier-Mukai transform for $H^*(X,{\Bbb Q})$),
there is an element $E_1$ of $M_H(w)$ such that
$\Hom(E_1,E) \ne 0$.
By induction hypothesis,
$E/E_1 \cong F_1 \oplus F_2 \oplus
\dots \oplus F_k$,
where $F_1$ is $S$-equivalent to $E_1^{\oplus n_1}$, $n_1 \geq 0$
and $F_i$, $i>1$ are $S$-equivalent to $E_i^{\oplus n_i}$.
Since $\Ext^1(E_i,E_1)=0$ for $i>1$,
$\Ext^1(\oplus_{i>1}F_i,E_1)=0$. 
Therefore $\oplus_{i>1}F_i$ is a direct summand of $E$,
which implies our claim.

Let $E_1$ be a stable sheaf of $v(E_1)=w$.
We set
\begin{equation}
 {\cal J}(l,E_1):=
 \{E \in {\cal M}(lw)^{ss}|\text{$E$ is $S$-equivalent to $E_1^{\oplus l}$}\}. 
\end{equation}
This is a closed substack of ${\cal M}(lw)^{ss}$
(see Appendix 4.3).
We shall next prove that
\begin{equation}\label{eq:j}
 \dim {\cal J}(l,E_1) \leq -1.
\end{equation}
Since $E_1$ is parametrized by the surface $M_H(w)$, \eqref{eq:j} implies that
\begin{equation}\label{eq:4.3}
 \dim {\cal M}(lw)^{ss}=l.
\end{equation}
For more details, see Appendix 4.3.
Proof of \eqref{eq:j}:
We set 
\begin{equation}
 {\cal J}(l,E_1,n):=\{E \in {\cal J}(l,E_1)|
 \dim \Hom(E_1,E)=n \}.
\end{equation}
By upper semi-continuity of cohomologies,
this is a locally closed substack of ${\cal J}(l,E_1)$.
If $n=l$, then
${\cal J}(l,E_1,l)=\{E_1^{\oplus l}\}$ and
it is a closed substack of ${\cal J}(l,E_1)$.
For an element $E$ of ${\cal J}(l,E_1,n)$,
there is an exact sequence
\begin{equation}
 0 \to \Hom(E_1,E) \otimes E_1 \to E \to E' \to 0
\end{equation}
where $E' \in {\cal J}(l-n,E_1,n')$.
The moduli number of $E$ which fits in this type of extension
is equal to $\dim{\cal J}(l-n,E_1,n')+nn'-n^2$.
Indeed, $\Hom(E_1,E) \otimes E_1 ( \cong E_1^{\oplus n})$ belongs to
${\cal J}(n,E_1,n)$ and $\dim{\cal J}(n,E_1,n)=-n^2$.
Hence by the proof of Lemma \ref{lem:A-est1}, we get the equality.
Therefore
we see that
\begin{equation}
 \dim {\cal J}(l,E_1) = -\min
 \begin{Sb}
  l=n_1+n_2+\dots+n_s\\
  n_1,n_2,\dots,n_s \geq 1
 \end{Sb}
 \left(\sum_{i=1}^s n_i^2-\sum_{i=1}^{s-1}n_in_{i+1} \right) \leq -1.
\end{equation}

\end{pf}
\begin{rem}
If $M_H(w)$ has a universal family, then Fourier-Mukai transform is defined
[Br].
Then ${\cal M}(lw)^{ss}$ is transformed to the stack of 0-dimensional
sheaves on $M_H(w)$.
In this case, by using [Y1, Thm. 0.4], we can get our lemma.
\end{rem}

\subsection{Lemma \ref{lem:-D*R} and its extensions}

\subsubsection{Proof of Lemma \ref{lem:-D*R}}
Proof of (1): By our assumptions, we have
\begin{equation}
 \frac{d_1}{r_1}=\frac{\deg E_1}{\rk E_1}<
 \frac{\deg \varphi(E_1)}{\rk \varphi(E_1)} \leq \frac{\deg E}{\rk E}=
 \frac{d}{r}.
\end{equation}
By Lemma \ref{lem:inequality}, $\deg \varphi(E_1)/\rk \varphi(E_1)=d/r$.
Hence $\varphi(E_1)$ coincides with $E$ except finite points of $X$.
Thus $\varphi$ is surjective in codimension 1.
We shall next prove that $\ker \varphi$ is a $\mu$-stable vector bundle.
Since $r_1-lr$ and $d_1-ld$ are relatively prime, we shall prove that
$\ker \varphi$ is $\mu$-semi-stable.
If $\ker \varphi $ is not $\mu$-semi-stable, then
there is a subsheaf $I$ of $\ker \varphi$
such that $I$ is $\mu$-stable and $d_1/r_1>\deg(I)/\rk(I)> (d_1-ld)/(r_1-lr)$.
Since $d/r>d_1/r_1$, we see that
\begin{equation}
 \frac{1}{r(r_1-lr)}>\frac{d}{r}-\frac{\deg(I)}{\rk (I)}
\geq \frac{1}{r \rk (I)}.
\end{equation}
Hence $\rk(I)>r_1-lr$, which is a contradiction.
Therefore $\ker \varphi$ is a $\mu$-stable vector bundle.

Proof of (2):
If $\phi$ is not injective,
then $(d_1-ld)/(r_1-lr)< \deg(\phi(E'))/\rk(\phi(E'))<d_1/r_1$.
Since $d/r>d_1/r_1$, we see that
\begin{equation}
 \frac{1}{r(r_1-lr)}>\frac{d}{r}-\frac{\deg(\phi(E'))}{\rk (\phi(E'))}
 \geq \frac{1}{r \rk (\phi(E'))}.
\end{equation}
Hence $\rk(\phi(E'))>r_1-lr$, which is a contradiction. 
Thus $\phi$ is injective.
Assume that $E$ is not $\mu$-semi-stable.
Then there is a $\mu$-stable quotient sheaf $F$ of $E$ 
such that $\deg F/\rk F<\deg E/\rk E$.
Since $F$ is also a quotient sheaf of $E_1$,
we have 
\begin{equation}
 \frac{d_1}{r_1}=\frac{\deg E_1}{\rk E_1}<\frac{ \deg F}{\rk F}
 <\frac{\deg E}{\rk E}=\frac{d}{r}.
\end{equation}
By Lemma \ref{lem:inequality}, we get
$\rk F \geq r_1+r$, which is a contradiction.
Therefore $E$ is $\mu$-semi-stable.
\qed

\begin{rem}
In order to define $\deg(E)$ in Lemma \ref{lem:-D*R}, 
we assumed that $\NS(X)={\Bbb Z}H$.
However $\deg(E)$ is still defined if
$(c_1(E),H)|(H,D)$ for all $D \in \NS(X)$.
Hence Lemma \ref{lem:-D*R} also holds under this assumption.
\end{rem}

\subsubsection{Extensions of Lemma \ref{lem:-D*R} and
[Y5, Lem. 4.1]}
We shall extend Lemma \ref{lem:-D*R} and
[Y5, Lem. 4.1].
Let $(X,H)$ be a polarized smooth projective surface of $\NS(X)={\Bbb Z}H$.
   
\begin{lem}\label{lem:codim=1II}
Let $E$ be a $\mu$-semi-stable vector bundle 
of $\rk(E)=lr$ and $\deg(E)=ld$
which is defined 
by a non-trivial extension 
\begin{equation}\label{eq:ext1}
 0 \to F_1 \to E \to F_2 \to 0,
\end{equation}
where $F_1$ and $F_2$ are $\mu$-stable vector bundles of
$\deg(F_1)/\rk(F_1)=\deg(F_2)/\rk(F_2)=d/r$.
Let $E_1$ be a $\mu$-stable vector bundle in
Lemma \ref{lem:-D*R}.
Let $\varphi: E_1 \to E$ be a non-trivial homomorphism.
Then $\varphi$ is surjective in codimension 1, or
$\varphi(E_1)$ is a subsheaf of $F_1$.
In particular, if $\Hom(E_1,F_1)=0$, then
$\ker \varphi$ is $\mu$-stable.
\end{lem}
\begin{pf}
We assume that $\varphi$ is not surjective in codimension 1.
By Lemma \ref{lem:inequality}, 
$\deg(\varphi(E_1))/\rk(\varphi(E_1))=d/r$.
Assume that $\varphi(E_1) \to F_2$ is not 0. Then
the $\mu$-stability of $F_2$ implies that it is surjective in codimension 1.
By the $\mu$-stability of $F_1$,
$F_1 \cap \varphi(E_1)=0$.
Thus we can regard $\varphi(E_1)$ as a subsheaf of $F_2$. 
Let $e \in \Ext^1(F_2,F_1)$ be the extension class of \eqref{eq:ext1}.
By the homomorphism $\Ext^1(F_2,F_1) \to \Ext^1(\varphi(E_1),F_1)$,
$e$ goes to $0$.
Since $F_1$ is a vector bundle,
$\Ext^1(F_2/\varphi(E_1),F_1)=0$.
Hence $\Ext^1(F_2,F_1) \to \Ext^1(\varphi(E_1),F_1)$ is injective.
Thus we get that $e=0$, which is a contradiction.
Hence $\varphi(E_1) \to F_2$ is a 0-map,
which means that $\varphi(E_1) \subset F_1$.
The last assertion follows from the proof of Lemma \ref{lem:-D*R}.
\end{pf}

\begin{lem}\label{lem:codim=1III}
Keep the notations in Lemma \ref{lem:codim=1II}. 
Assume that $\Ext^1(E_1,F_2)=0$.
Then a non-trivial extension of $E$ by $E_1$ is $\mu$-stable.
\end{lem}
\begin{pf}
Let $E'$ be a non-trivial extension
of $E$ by $E_1$.
\begin{equation}\label{eq:E'1}
0 \to E_1 \to E' \to E \to 0.
\end{equation}
We shall prove that $E'$ is $\mu$-stable.
We consider the following diagram which is induced by
the extension \eqref{eq:E'1}.
\begin{equation}\label{eq:E'2}
 \begin{CD}
  @. @. 0 @. 0 @.\\
  @.@. @AAA @AAA @.\\
  @. @. F_2 @= F_2 @.\\
  @.@. @AAA @AAA @.\\
  0 @>>> E_1 @>>> E' @>>> E @>>> 0 \\
  @. @| @AAA @AAA @.\\
  0 @>>> E_1 @>>> E'' @>>> F_1 @>>> 0 \\
  @.@. @AAA @AAA @.\\
  @. @. 0 @. 0 @.
 \end{CD}
\end{equation}
By Serre duality, $\Ext^1(F_2,E_1)=0$.
Hence $\Ext^1(E,E_1) \to \Ext^1(F_1,E_1)$ is injective.
Thus the last horizontal sequence of \eqref{eq:E'2} does not split.
By Lemma \ref{lem:4-1}, $E''$ is $\mu$-stable.
If the middle vertical sequence splits, then
$E \cong F_1 \oplus F_2$, which is a contradiction.
Hence $E'$ is a non-trivial extension of $F_2$ by a
locally free sheaf $E''$.
By the construction of $E''$, the conditions in
Lemma \ref{lem:4-1} hold.
Hence applying Lemma \ref{lem:4-1} again, we see that
$E'$ is $\mu$-stable.
\end{pf} 

\subsubsection{Extension of [Y5, Lem. 4.2]}
Let $(X,H)$ be a polarized K3 surface
of $\Pic(X)={\Bbb Z}H$.
Let $E_1$ be an exceptional vector bundle of
$v(E_1):=r_1+d_1 H+a_1 \omega$
and let $v=l(r+dH)+a \omega$ be a primitive
Mukai vector of $d r_1-d_1r=1$.  
We set 
\begin{equation}
 M_H(v)_i^{\mu s}:=\{ E \in M_H(v)^{\mu s}|
 \dim \Hom(E_1,E)=i-\langle v,v(E_1) \rangle \}.
\end{equation}
Assume that $r_1>lr$.
Then we get the following estimate which is necessary for the condition 
section 0.2 (a).
\begin{lem}\label{lem:codim=1I}
\begin{enumerate}
\item[(1)]
If $\langle v,v(E_1) \rangle <0$,
then 
\begin{equation}
 \codim_{M_H(v)^{\mu s}} M_H(v)_i^{\mu s}  \geq -\langle v,v(E_1) \rangle+1
 \geq 2
\end{equation}
for $i \geq 1$.
\item[(2)]
If $\langle v,v(E_1) \rangle  \geq 0$,
then 
\begin{equation}
 \codim_{M_H(v)^{\mu s, loc}} (M_H(v)_i^{\mu s}\cap
 M_H(v)^{\mu s, loc}) \geq \langle v,v(E_1) \rangle+1
 \geq 1
\end{equation}
for $i>\langle v,v(E_1) \rangle$,
where $M_H(v)^{\mu s, loc}=M_H(v)^{\mu s}\cap M_H(v)^{loc}$.
\end{enumerate}
\end{lem}
\begin{pf}
We set 
\begin{equation}
N_i:= \left\{E_1 \subset E\left|
\begin{aligned}
&\text{$E \in M_H(u)$, $\dim \Hom(E_1,E)=i+1-\langle v,v(E_1) \rangle $ },\\
&\text{$E/E_1$ is a $\mu$-stable sheaf of $v(E/E_1)=v$}
\end{aligned}
\right.
\right \},
\end{equation}
where $u=v+v(E_1)$.
Then we see that $\dim N_i \leq
i-\langle v, v(E_1) \rangle+\dim M_H(u)=i+\langle v,v(E_1) \rangle
+\dim M_H(v)-2$.
Let $\pi_v': N_i  \to M_H(v)^{\mu s}$ be the
morphism
sending $(E_1 \subset  E) \in  N_i$ to
$E/E_1 \in  M_H(v)^{\mu s}$.
By Lemma  \ref{lem:4-1}, $\pi_v'(N_i)=
M_H(v)_i^{\mu s}$ and 
${\pi_v'}^{-1}(E/E_1)$ is isomorphic
to
the projective space 
${\Bbb P}(\Ext^1(E/E_1,E_1)^{\vee})$.
Hence we get that
$\dim M_H(v)_i^{\mu s}=\dim N_i-(i-1) \leq
\dim M_H(v)+\langle v,v(E_1) \rangle-1$.
Thus (1) holds.

We next prove the second claim.
Hence we assume that 
$i>\langle v,v(E_1) \rangle$.
For $E \in M_H(v)_i^{\mu s}\cap
M_H(v)^{\mu s, loc}$,
we choose a homomorphism $\phi:E_1 \to E$.
By Lemma \ref{lem:-D*R}, $\phi$ is surjective in codimension 1.
We set $G:=\coker(E^{\vee} \to E_1^{\vee})$.
By Lemma \ref{lem:-D*R},
 $F$ is a stable sheaf of $v(G)=w:=v(E_1)^{\vee}-v^{\vee}$.
It is easy to see that $\dim \Hom(E_1^{\vee},G)
=\dim \Ext^1(E_1,E)+1=\langle v,v(E_1) \rangle+i+1$.  
Hence $\phi:E_1 \to E$ is parametrized by an open subscheme of a 
projective bundle of dimension $(\langle v,v(E_1) \rangle+i)$
over the subscheme $M_H(w)_i$,
where
\begin{equation}
 \begin{split}
  M_H(w)_i:&=\{G \in M_H(w)|\dim \Hom(E_1^{\vee },G)=
  \langle v,v(E_1) \rangle+i+1\}\\
  &=\{G \in M_H(w)|\dim \Hom(E_1^{\vee },G)=\langle w,v(E_1)^{\vee} \rangle
  +(i-1)\}.
 \end{split}
\end{equation}
Thus we see that
\begin{equation}
 \begin{split}
  \dim M_H(v)_i^{\mu s}\cap
  M_H(v)^{\mu s, loc}& \leq \langle w^2 \rangle+2
  +(\langle v,v(E_1) \rangle+i)-(i-1)\\
  &=\langle v^2 \rangle+2-(\langle v,v(E_1)\rangle+1).
 \end{split}
\end{equation}
 \end{pf}

\section{Case A}

\subsection{Estimate}

In this section, we fix a primitive Mukai vector 
$r+\xi, \xi \in \NS(X)$. 
We assume that 
\begin{equation}\label{eq:A}
((\xi^2)+2)/2r \not \in {\Bbb Z}.
\end{equation}
We shall prove Theorem \ref{thm:main} for a primitive
Mukai vector $v:=l(r+\xi)+a\omega \in H^*(X,{\Bbb Z})$. 
We shall first estimate the dimensions of various locally closed substacks of
${\cal M}(v)$.

\begin{lem}\label{lem:non-empty}
If ${\cal M}(v)^{\mu ss} \ne \emptyset$, then 
$\langle v^2 \rangle \geq 0$.
If the equality holds, then ${\cal M}(v)^{\mu ss}={\cal M}(v)^{ss}$.
\end{lem}
\begin{pf}
Let $E$ be a $\mu$-semi-stable sheaf of $v(E)=v$
and $E$ is $S$-equivalent to $\oplus_{i=1}^s E_i$
with respect to $\mu$-stability,
where $E_i, 1 \leq i \leq s$ are $\mu$-stable sheaves.
We set
\begin{equation}
v(E_i):=l_i(r+\xi)+a_i\omega, 1 \leq i \leq s.
\end{equation}
By our assumption \eqref{eq:A},
$\langle v(E_i)^2 \rangle=l_i(l_i(\xi^2)-2r a_i) \ne -2$.
Thus $\langle v(E_i)^2 \rangle \geq 0$ for all $i$.
Since
\begin{equation}
\frac{\langle v^2 \rangle}{l}=
\sum_{i=1}^s\frac{\langle v(E_i)^2 \rangle}{l_i},
\end{equation}
we get $\langle v^2 \rangle \geq 0$.
If $\langle v^2 \rangle=0$, then 
$\langle v(E_i)^2 \rangle =0$ for all $i$.
Since $\langle v(E_i)^2 \rangle/\rk(E_i)^2=(\xi^2)-2 a_i/r l_i$ and
$\chi(E_i)/\rk(E_i)=1+a_i/rl_i$,
we see that $\chi(E_i)/\rk(E_i)=\chi(E)/\rk(E)$ for all $i$.
Thus $E$ is semi-stable.
\end{pf}
\begin{cor}\label{cor:isotropic}
If $\langle v^2 \rangle=0$, then ${\cal M}(v)^{\mu ss}$ 
consists of locally free sheaves.
\end{cor}

\begin{defn}
Let $w=l_0(r+\xi)+a_0 \omega$ be the primitive Mukai vector  
such that $\langle w^2 \rangle =0$.
\end{defn}
By Lemma \ref{lem:non-empty} and Corollary \ref{cor:isotropic},
 $M_H(w)$ consists of $\mu$-stable locally free sheaves.

\begin{lem}\label{lem:stable}
\begin{enumerate}
\item[(1)]
\begin{equation}
\dim ({\cal M}(v)^{\mu ss} \setminus {\cal M}(v)^{ss})
\leq \langle v^2 \rangle.
\end{equation}
\item [(2)]
Assume that $\langle v^2 \rangle>0$. Then
\begin{equation}
\dim ({\cal M}(v)^{\mu ss} \setminus {\cal M}(v)^{s})
\leq \langle v^2 \rangle.
\end{equation}
In particular,
if ${\cal M}(v)^{\mu ss} \ne \emptyset$, then
${\cal M}(v)^s \ne \emptyset$ and
$\dim {\cal M}(v)^{\mu ss}=\langle v^2 \rangle+1$.
\end{enumerate}
\end{lem}
\begin{pf}
By Lemma \ref{lem:pss}, it is sufficient to prove (1).
Let $F$ be a $\mu$-semi-stable sheaf of $v(F)=v$.
We assume that $F$ is not stable.
Let 
\begin{equation}
 0 \subset F_1 \subset F_2 \subset \dots \subset F_s=F
\end{equation}
be the Harder-Narasimhan filtration of $F$.
We set 
\begin{equation}
 v_i:=v(F_i/F_{i-1})=l_i(r+\xi)+a_i \omega, 1 \leq i \leq s.
\end{equation}
Since $\chi(F_i/F_{i-1})/\rk(F_i/F_{i-1})>
\chi(F_{i+1}/F_{i})/\rk(F_{i+1}/F_{i})$,
we get that
\begin{equation}\label{eq:a*}
 \frac{a_0}{l_0} \geq \frac{a_1}{l_1}>\frac{a_2}{l_2}> \dots
 >\frac{a_s}{l_s}.
\end{equation}
Let ${\cal F}^{HN}(v_1,v_2,\dots,v_s)$ be the substack of ${\cal M}(v)^{\mu ss}$
whose element $E$ has the Harder-Narasimhan filtration
of the above type.
We shall prove that
$\dim {\cal F}^{HN}(v_1,v_2,\dots,v_s) \leq \langle v^2 \rangle$.
Since $\Hom(F_i/F_{i-1},F_j/F_{j-1})=0$ for $i<j$,
Lemma \ref{lem:A-est2} in Appendix implies that
\begin{equation}
  \dim {\cal F}^{HN}(v_1,v_2,\dots,v_s)=
  \sum_{i=1}^s \dim {\cal M}(v_i)^{ss}+
  \sum_{i<j}\langle v_j,v_i \rangle.
\end{equation}
For $i<j$, by using Lemma \ref{lem:non-empty} and \eqref{eq:a*},
we see that
\begin{equation}\label{eq:v_i,v_j}
\begin{split}
\langle v_i,v_j \rangle &=
l_il_j (\xi^2)-(l_ia_j+l_j a_i)r\\
&=l_il_j (\xi^2)-2l_ja_i r+
(a_il_j-a_jl_i)r\\
&=l_j(l_i (\xi^2)-2a_ir)+(a_il_j-a_jl_i)r\\
& \geq (a_il_j-a_jl_i)r \geq r \geq 2,
\end{split}
\end{equation}
where the inequality $r \geq 2$ comes from our assumption \eqref{eq:A}.
Hence if $\langle v_i^2 \rangle>0$ for all $i$,
then, by using Lemma \ref{lem:pss},
we see that
\begin{equation}
 \dim {\cal F}^{HN}(v_1,v_2,\dots,v_s) \leq
 (\langle v^2 \rangle+1)-
 \left(\sum_{i<j}\langle v_i,v_j \rangle-s+1 \right)
 \leq \langle v^2 \rangle.
\end{equation} 
Assume that $\langle v_i^2 \rangle =0$, i.e. $v_i=l_i'w$, $l_i \in {\Bbb Z}$.
Then $i=1$ and $(a_1l_j-a_jl_1)$ is divisible by $l_i'$.
Hence 
\begin{equation}\label{eq:v_1,v_j}
 \begin{split}
  \langle v_1,v_j \rangle-l_1' &=l_1'(\langle w,v_j \rangle-1)\\
  & \geq l_1'(r-1)>0.
 \end{split}
\end{equation}
In this case, by using Lemma \ref{lem:ss}, we see that
\begin{equation}
 \dim {\cal F}^{HN}(v_1,v_2,\dots,v_s) \leq
 (\langle v^2 \rangle+1)-
 \left(\sum_{i<j}\langle v_i,v_j \rangle-(l_1'+s-1)+1\right)
 \leq \langle v^2 \rangle.
\end{equation} 
Hence we get our lemma.
\end{pf}

\begin{prop}\label{prop:stable}
Assume that $\langle v^2 \rangle >0$. Then
\begin{equation}
\dim( {\cal M}(v)^s \setminus {\cal M}(v)^{\mu s})
\leq \langle v^2 \rangle.
\end{equation}
In particular,
${\cal M}(v)^{\mu s} \ne \emptyset$, if ${\cal M}(v)^{s} \ne \emptyset$.
\end{prop}
\begin{pf}
Let $E$ be a stable sheaf and $E_1$ be a $\mu$-stable subsheaf of $E$ such that
$E/E_1$ is torsion free.
We set
\begin{equation}
 \begin{split}
  v_1:=v(E_1)=l_1(r+\xi)+a_1 \omega,\\
  v_2:=v(E/E_1)=l_2(r+\xi)+a_2 \omega.
 \end{split}
\end{equation}
Since $\chi(E_1)/\rk E_1<\chi(E)/\rk E$,
we get $\langle v(E_1)^2 \rangle>0$ and 
\begin{equation}\label{eq:stable}
 \frac{a_1}{l_1}<\frac{a_2}{l_2}.
\end{equation}
Let $J(v_1,v_2)$ be the substack  of
${\cal M}(v)^s$ consisting of $E$ which has a subsheaf
$F_1 \subset E$.
By using Lemma \ref{lem:A-est1} in Appendix, we shall estimate
$\dim J(v_1,v_2)$.
By Lemma \ref{lem:hom},
$ \dim \Hom(E_1,E/E_1) \leq l_2/l_1$,
 and if $E_1$ is general, then $\Hom(E_1,E/E_1)=0$.
We shall bound the dimension of the substack
\begin{equation}
 {\cal N}(v_1,v_2):=\{(E_1, E_2) \in {\cal M}(v_1)^{\mu ss} \times
 {\cal M}(v_2)^{\mu ss}|\dim \Hom(E_1,E_2) \ne 0 \}.
\end{equation}
For a fixed $E_2 \in {\cal M}(v_2)^{ss}$,
\begin{equation}
 \#\{E_1^{\vee \vee}|
 E_1 \in {\cal M}(v_1)^{\mu s}, \Hom(E_1,E_2) \ne 0 \} <\infty.
\end{equation}
Hence, by using [Y1, Thm. 0.4],
we see that
\begin{equation}
 \dim \{E_1 \in {\cal M}(v_1)^{\mu s}|\Hom(E_1,E_2) \ne 0 \}
 \leq \dim {\cal M}(v_1)^{\mu s}-2-(\rk v_1-1).
\end{equation}
Thus $\dim {\cal N}(v_1,v_2) \leq \dim{\cal M}(v_1)^{\mu ss}+
\dim {\cal M}(v_2)^{\mu ss}-3$.
Moreover, taking \eqref{eq:JHF} into account,
 if $l_1 \ne l_0$ and $v_2=l_2'w, 
l_2' \in {\Bbb Z}$, then we get
${\cal N}(v_1,v_2)=\emptyset$.

If $\langle v_2^2 \rangle >0$, then
Lemma \ref{lem:stable} implies that
$\dim {\cal M}(v_2)^{\mu ss}=\langle v_2^2 \rangle+1$.
Hence Lemma \ref{lem:A-est1} implies that 
\begin{equation}
 \begin{split}
   \dim {\cal M}(v)^s-\dim J(v_1,v_2)&=
   \min \left(\langle v_1,v_2\rangle-\frac{l_2}{l_1}+2,
   \langle v_1,v_2\rangle-1 \right)\\
   &=l_1 \frac{\langle v_2^2 \rangle}{2l_2}+
   l_2\frac{\langle v_1^2 \rangle}{2l_1}-
   \max\left(\frac{l_2}{l_1}-2,1 \right)>0.
 \end{split}
\end{equation}
We next treat the case where $\langle v_2^2 \rangle=0$.
Then $v_2=l_2' w, l_2' \in {\Bbb Z}$.
By Lemma \ref{lem:stable} (1) and Lemma \ref{lem:ss},
$\dim {\cal M}(v_2)^{\mu ss}=\langle v_2^2 \rangle+l_2'$.
If $l_1=l_0$, then $l_2/l_1=l_2'$. So we see that
\begin{equation}
 \begin{split}
   \dim {\cal M}(v)^s-\dim J(v_1,v_2)&=
   \min \left(l_2'(\langle v_1,w\rangle-1-1)+3,
   l_2'(\langle v_1,w\rangle-1) \right)\\
   &=\min \{l_2'((a_0l_1-a_1l_0)r-2)+3,l_2'((a_0l_1-a_1l_0)r-1)\} >0.
 \end{split}
\end{equation}
\end{pf}
\begin{rem}\label{rem:codim=1}
By the proof of Lemma \ref{lem:stable} and Proposition \ref{prop:stable},
we see that 
\begin{equation}
 \codim_{{\cal M}(v)^{\mu ss}}({\cal M}(v)^{\mu ss}
 \setminus {\cal M}(v)^{\mu s}) \geq 2 \text{ for $r \geq 3$}.
\end{equation}
Moreover if $r=2$, then 
the general member $E$ of ${\cal M}(v)^{s}
\setminus {\cal M}(v)^{\mu s}$ fits in the following exact sequence
\begin{equation}\label{eq:codim=1}
 0 \to E_1 \to E \to E_2 \to 0
\end{equation}
where $E_1$ is a $\mu$-stable vector bundle
and $E_2$ is a $\mu$-stable vector bundle of $v(E_2)=w$. 

Indeed, if $\codim_{{\cal M}(v)^{s}}({\cal M}(v)^{s}
\setminus {\cal M}(v)^{\mu s})=1$, then
(1) $\langle v_2^2 \rangle>0$, $l_1=l_2=\langle v_1^2 \rangle/2=
\langle v_2^2 \rangle/2=1$, or
(2) $r=2$, $l_2'=1$ and $a_0l_1-a_1l_0=1$.
By \eqref{eq:stable}, case (1) does not occur.  
\end{rem}

\subsection{Proof of Theorem \ref{thm:main} for case A}

\subsubsection{The case of $r>2$}

In the same way as in the proof of Theorem \ref{thm:period},
we shall prove Theorem \ref{thm:main}
if $r >2$.
Let $v=l(r+\xi)+a \omega$ be a primitive Mukai vector 
on a K3 surface $X$ such that $r+\xi, \xi \in  \NS(X)$ is primitive.
We claim that we can find  
a primitive Mukai vector $v'=l(r+d'H')+a' \omega$ on a polarized 
K3 surface $(X',H')$ of $\Pic(X')={\Bbb Z}H'$ 
and an exceptional vector bundle $G$ of $v(G):=v_1=
r_1+d_1H'+a_1 \omega'$
such that
(1) $r+d' H'$ is primitive,
(2) $\langle v^2 \rangle=\langle (v')^2 \rangle$,
(3) $a \equiv a' \mod l$,
(4) $d'r_1-d_1r=1$,
(5) $r_1-lr \geq 2$ and 
(6) $\langle v(G),v' \rangle=-1$,
where $\omega'$ is the fundamental class of $X'$.
\newline
Proof of the claim:
We choose integers $r_1,d_1$ and $d'$
such that $ar_1 \equiv 1 \mod l$ and 
$(r_1,r)=1$.
We can easily choose such an integer $r_1$ of 
$r_1-lr \geq 2$.
Then we can choose $d'$ and $d_1$ of $d' r_1-d_1r=1$.
Let $q$ be an integer of $a r_1+ql=1$.
We set $k(s):=r_1(qr+r_1s)-r^2, s=(\xi^2)/2 \in {\Bbb Z}$.
Then $k(s)=(r_1 \langle v^2 \rangle/2l+r)r_1/l-r^2$.
Since $\langle v^2 \rangle \geq 0$, $k(s) \geq rr_1^2/l-r^2>0$.
Let $(X',H')$ be a polarized K3 surface such that 
$\Pic(X')={\Bbb Z} H'$ and $({H'}^2)=2k(s)$.
We set
\begin{equation}\label{eq:4-1}
\begin{cases}
v'=lr+ld' H'+\{l((1+d'r_1)d_1s+{d'}^2qr_1-r{d'}^2)+a \}\omega',\\
v_1=r_1+d_1 H'+\{r_1(-{d'}^2+d_1^2s)+d_1^2rq+2d' \}\omega'.
\end{cases}
\end{equation}
Then we see that
\begin{equation}
\begin{cases}
\langle v_1^2 \rangle=-2,\\
\langle {v'}^2 \rangle =2l(ls-ra)=\langle {v}^2 \rangle,\\
\langle v_1,v' \rangle=-1.
\end{cases}
\end{equation}
Since $r_1$ and $d_1$ are relatively prime,
Theorem \ref{thm:period} implies that there is an exceptional vector bundle
$G$ of $v(G)=v_1$.   
Then $G$ and $v'$ satisfy our claims.

We shall consider reflection defined by $G$. 
We note that
\begin{equation}
w:=-R_{v_1}(v')^{\vee}=(r_1-lr)-(d_1-ld')H'+(a_1-a')\omega'.
\end{equation}
Since $r_1-lr$ and $d_1-ld'$ are relatively prime,
Theorem \ref{thm:period} for a Mukai vector $w$ of $\ell(w)=1$
implies that $M_{H'}(w) \ne \emptyset$
and Theorem \ref{thm:main} holds for this space.
Since $r_1-lr \geq 2$ and $M_{H'}(w)$ consists of $\mu$-stable sheaves,
[Y1, Thm. 0.4] implies that there is a $\mu$-stable vector bundle 
$E$ of $v(E)=w$.
By (6), we see that $\chi(G^{\vee},E)=-\langle v(G)^{\vee},w \rangle=1$.
Hence there is a non-trivial homomorphism $\phi:G^{\vee} \to E$.
By using Lemma \ref{lem:-D*R} (2), we see that
$\coker( \phi^{\vee})$ is a $\mu$-semi-stable sheaf
of $v(\coker( \phi^{\vee}))=v'$.
Then Lemma \ref{lem:stable} and 
Proposition \ref{prop:stable} imply that 
$M_{H'}(v')^{\mu s} \ne \emptyset$.
Applying Lemma \ref{lem:-D*R} and
Lemma \ref{lem:codim=1I},
we get a birational map
\begin{equation}\label{eq:birat}
M_{H'}(v') \cdots \to M_{H'}(w)
\end{equation}
sending $F \in M_{H'}(v')_0^{\mu s} \cap M_{H'}(v')^{loc}$
to $\coker(F^{\vee} \to G^{\vee})$,  
which means that Theorem \ref{thm:main} (1) and (2-1) hold for
$M_{H'}(v')$.
Then Proposition \ref{prop:deform} implies that
Theorem \ref{thm:main} (1) and (2-1) also hold for $M_H(v)$,
where $H$ is a general ample divisor on $X$.

Moreover by Remark \ref{rem:codim=1},
we can naturally identify $H^2(M_H(v'),{\Bbb Z})$ with
$H^2(M_H(w),{\Bbb Z})$.
Let ${\cal F}^H:H^*(X,{\Bbb Z}) \to H^*(X,{\Bbb Z})$ be an isomerty of
Mukai lattice defined by
${\cal F}^H(x)=R_{v_1}(x)^{\vee}$.
Then it is easy to see that the following diagram is commutative
(see the computation in [Y5, 2.4]).
\begin{equation*}
\begin{CD}
{v'}^{\perp} @>{{\cal F}^H}>> w^{\perp}\\
@V{\theta_{v'}}VV  @VV{\theta_{w}}V\\
H^2(M_{H'}(v'),{\Bbb Z}) @=
H^2(M_{H'}(w),{\Bbb Z})
\end{CD}
\end{equation*}
Hence (2-2) also holds if $r>2$.
\qed

\subsubsection{The case of $r=2$}

We next treat the case where $r=2$.
It is sufficient to extend the birational map
\eqref{eq:birat} to a general member $E \in M_{H'}(v')$ which fits in
\eqref{eq:codim=1}. 
Let $G$ be the exceptional vector bundle in 2.2.1.
By using Lemma \ref{lem:codim=1II} and  \ref{lem:codim=1III},
we can prove the following:
Assume that $\langle v',v(G) \rangle=-1$.
Then, for a general member $E$ which fits in the exact sequence
\eqref{eq:codim=1},
\begin{enumerate}
\item[(1)]
$\Ext^1(G,E)=0$,
\item[(2)]
$\phi:\Hom(G,E) \otimes G \to E$ is surjective in codimension 1 and
$\ker \phi$ is stable.
\end{enumerate}
\begin{pf}
Since $\langle (v')^2 \rangle >0$, we can write 
$v'=x w -y\omega$, where $x,y  \in {\Bbb Q}$ and $x, y >0$.
Since $\langle v',v(G) \rangle=x \langle w, v(G) \rangle+y \rk(G)=-1$,
$\langle w,v(G) \rangle<0$.
Hence $\langle v(E_1), v(G) \rangle \geq 0$ and
$\langle v(E_2), v(G) \rangle < 0$.
Applying Lemma \ref{lem:codim=1I},
we may assume that $\Hom(G,E_1)=\Ext^1(G,E_2)=0$.
By Lemma \ref{lem:codim=1III},
we can use the same argument as in the proof of 
Lemma \ref{lem:codim=1I} (1).
Hence (1) holds.
Applying Lemma \ref{lem:codim=1II}, we get (2).
\end{pf}
Therefore the dual of $\phi$ : 
$E^{\vee} \to G^{\vee}$ is injective and the cokernel is stable.
Thus the reflection induces a birational map $M_H(v) \setminus Z
 \to M_H(v(G)^{\vee}-v^{\vee})$ such that
$\codim_{M_H(v)} Z \geq 2$.
Therefore Theorem \ref{thm:main} also holds for this case.
\qed

\begin{rem}\label{rem:A}
Under the assumptions for case A, the following holds.
\begin{enumerate}
\item [(1)]
$M_H(v) \ne \emptyset$ if and only if $\langle v^2 \rangle \geq 0$.
\item [(2)]
$M_H(v)^{\mu s} \ne \emptyset$ if and only if $\langle v^2 \rangle \geq 0$.
\end{enumerate}
In particular, $M_H(v)$ contains $\mu$-stable vector bundles
if $M_H(v) \ne \emptyset$.
\end{rem}

\section{Case B}

In this section, we assume that there is a $(-2)$ vector $v_0$ of the form
$v_0=r+\xi+b \omega, b \in {\Bbb Z}$.
Let $E_0$ be the element of $M_H(v_0)$.
We shall prove Theorem \ref{thm:main} for a primitive Mukai vector
$v:=l v_0-a \omega$.
If $\langle v^2 \rangle>2l^2$, then the same proof in section 2.2.1 works,
because of Lemma \ref{lem:quot}.
Hence 
we may assume that $\langle v^2 \rangle \leq 2l^2$. 

\subsection{The case of $\langle v^2 \rangle<2l^2$ }

We first treat the case where $\langle v^2\rangle<2l^2$.
Clearly $\langle v^2 \rangle \geq -2$, if $M_H(v) \ne \emptyset$.
If $\langle v^2 \rangle=-2$, then $2l(a \rk v_0-l)=-2$.
Hence we get $l=1$.
This case is covered in Theorem \ref{thm:period}.
So we assume that $\langle v^2 \rangle \geq 0$, that is
\begin{equation}\label{eq:small}
l \leq a \rk v_0 < 2l.
\end{equation}
Based on the next key lemma, we shall prove Theorem \ref{thm:main}
in 3.1.2 and 3.1.3. 
\begin{lem}\label{lem:not-free}
$\codim_{M_H(v)}(M_H(v) \setminus M_H(v)^{loc}) \geq 2$
unless 
\begin{equation}
 v=
  \begin{cases}
   (\rk v_0)v_0-\omega,\\
   lv_0-(l+1) \omega,\;\rk v_0=1.
  \end{cases}
\end{equation}
\end{lem}

\subsubsection{Proof of Lemma \ref{lem:not-free}}

We set
\begin{equation}
 {\cal M}(v,n)^s:=\{E \in {\cal M}(v)^s|\dim (E^{\vee \vee}/E)=n \}.
\end{equation}
For $E \in {\cal M}(v,n)^s$,
let 
\begin{equation}
0 \subset F_1 \subset F_2 \subset \dots \subset F_s=E^{\vee \vee}
\end{equation}
be the Harder-Narasimhan filtration of $E^{\vee \vee}$.
We set 
\begin{equation}
v_i:=v(F_i/F_{i-1})=l_iv_0-a_i \omega, 1 \leq i \leq s.
\end{equation}
We shall estimate the codimension of 
the substack
\begin{equation}
 {\cal M}(v,n;v_1,v_2,\dots,v_s):=
 \{E \in {\cal M}(v,n)^s|E^{\vee \vee} 
 \in {\cal F}^{HN}(v_1,v_2,\dots,v_s)\},
\end{equation}
where ${\cal F}^{HN}(v_1,v_2,\dots,v_s)$ is defined as in the proof of
Lemma \ref{lem:stable}. 
We divide our consideration into two cases
\begin{enumerate}
\item[(I)] $s \geq 2$, that is, $E^{\vee \vee}$ is not semi-stable.
\item[(II)] $s=1$, that is, $E^{\vee \vee}$ is semi-stable.
\end{enumerate}
Case (I).
(I-a) If $\langle v_1^2 \rangle \geq 0$, then 
we can use almost the same arguments as in the proof of Lemma \ref{lem:stable}.
The difference comes from the inequality $r \geq 2$
which was used in \eqref{eq:v_i,v_j} and \eqref{eq:v_1,v_j}.
Thus 
\begin{equation}
\dim {\cal F}^{HN}(v_1,v_2,\dots,v_s) \leq
\langle v(E^{\vee \vee})^2 \rangle+1.
\end{equation}
The equality holds only if $r=1$ and $\langle v_1^2 \rangle=0$.
Since $F_s$ is locally free, $F_1$ is also locally free.
On the other hand, if $r=1$ and $\langle v_1^2 \rangle=0$,
then ${\cal M}(v_1)^{ss}$ 
consists of non-locally free sheaves (cf. Lemma \ref{lem:ss}).
Therefore the equality does not hold.
Hence, by using [Y1, Thm. 0.4], we see that 
\begin{equation}
 \begin{split}
  \dim {\cal M}(v,n;v_1,v_2,\dots,v_s) & <
  \langle (v-n\omega)^2 \rangle+1+n(\rk v+1)\\
  & \leq \langle v^2 \rangle +1-n(\rk v-1).
 \end{split}
\end{equation}
Thus we get a desired estimate
\begin{equation}
\codim_{{\cal M}(v)^s}{\cal M}(v,n;v_1,v_2,\dots,v_s) \geq 2.
\end{equation}

(I-b) We assume that $\langle v_1^2 \rangle < 0$.
Then $F_1=E_0^{\oplus l_1}$, which implies that
$\dim {\cal M}(v_1)^{ss}=-l_1^2=\langle v_1^2 \rangle+l_1^2$.
For convenience sake,
we set 
\begin{equation}
 v_2':=\sum_{i=2}^s v_i=l_2'v_0-a_2' \omega.
\end{equation}
(I-b-1) We first assume that 
\begin{equation}\label{eq:v_2}
 \langle (v_2')^2 \rangle
 =2l_2'(a_2' \rk v_0-l_2')>0.
\end{equation}
Then, since $\langle v_i^2 \rangle \geq 0$ for all $i \geq 2$,
 $\dim {\cal F}^{HN}(v_2,\dots, v_s) \leq \langle  (v_2')^2 \rangle+1$.
Thus we get
\begin{equation}
\langle v(E^{\vee \vee})^2 \rangle+1-\dim {\cal F}^{HN}(v_1,v_2,\dots,v_s)
 \geq \langle v_1,v_2' \rangle-l_1^2.
\end{equation}
We shall prove that 
\begin{equation}\label{eq:codim}
\langle v_1,v_2' \rangle-l_1^2+n(l \rk v_0-1) \geq 2.
\end{equation}
\begin{pf}
Since $E$ is stable, $F_1 \cap E$ satisfies that
\begin{equation}
\frac{\chi(F_1 \cap E)}{\rk F_1}<\frac{\chi(E)}{\rk E}.
\end{equation}
Since $\chi(F_1 \cap E) \geq \chi(F_1)-n$ and $v_1=l_1 v_0$
(i.e. $a_1=0$), we see that
\begin{equation}\label{eq:n1}
nl_2'-a_2'l_1>0.
\end{equation}
By our assumption \eqref{eq:small},
$(a_2'+n) \rk v_0 < 2(l_1+l_2')$.
By using \eqref{eq:v_2}, we see that
\begin{equation}\label{eq:n2}
 n \leq \frac{2l_1+l_2'-1}{\rk v_0}.
\end{equation}
By using \eqref{eq:v_2} and \eqref{eq:n1}, we see that
\begin{equation}\label{eq:n3}
 nl_1 \rk v_0 \geq l_1^2+1.
\end{equation}
Assume that $\rk v_0 \geq 2$. Then \eqref{eq:n2} implies that
$n<l_1+l_2'/2$. By using \eqref{eq:v_2} and \eqref{eq:n3},
we see that
\begin{equation}
 \begin{split}
  \langle v_1,v_2' \rangle-l_1^2+n(l \rk v_0-1)=
  &l_1(-2l_2'+a_2' \rk v_0)+n((l_1+l_2')\rk v_0-1)-l_1^2\\
  =&-l_1l_2'+l_1(a_2' \rk v_0-l_2')+nl_1 \rk v_0
  +l_1l_2'+l_2'(n\rk v_0-l_1)-n-l_1^2\\
  \geq &l_1(a_2' \rk v_0-l_2')+l_2'(n \rk v_0-l_1)+1-n \\
  \geq & l_1+l_2'+1-n>1.
 \end{split}
\end{equation}
Assume that $\rk v_0=1$.
Then similar computations work if $a_2' \rk v_0-l_2'>1$.
So we assume that $\rk v_0=a_2' \rk v_0-l_2'=1$.
Then \eqref{eq:n3} implies that $n-l_1>0$.
Hence we see that
\begin{equation}\label{eq:l}
 l_1(-2l_2'+a_2' \rk v_0)+n((l_1+l_2')\rk v_0-1)-l_1^2
 = (l_1+l_2'-1)(n-l_1) \geq (l-1) \geq 1.
\end{equation}
If the equality holds, then $l=2$ and $n-l_1=1$.
In this case, we get that $l_1=l_2'=1$ and $n=2$.
By \eqref{eq:n1}, we get a contradiction.
Thus the left hand side of \eqref{eq:l} is greater than or equal to 2.
Therefore \eqref{eq:codim} holds.
\end{pf}
By \eqref{eq:codim} and [Y1, Thm. 0.4], we get a desired estimate 
\begin{equation}
 \codim_{{\cal M}(v)^s} {\cal M}(v,n;v_1,v_2,\dots,v_s) \geq 2.
\end{equation}

(I-b-2) We next assume that
\begin{equation}\label{eq:v_2-2}
 \langle (v_2')^2 \rangle
 =2l_2'(a_2' \rk v_0-l_2')=0.
\end{equation}
Then $s=2$ and $\langle v_0,v_2' \rangle=-l_2'$.
Since $\Hom(E_0, F_2/F_1)=0$, 
$\dim \Hom(F_2/F_1,E_0) \geq l_2'$.
By Lemma \ref{lem:hom}, coevaluation map
\begin{equation}
 \psi:F_2/F_1 \to E_0 \otimes \Hom(F_2/F_1,E_0)^{\vee}
\end{equation}
is surjective in codimension 1 (cf. Lemma \ref{lem:hom}).
Therefore we get 
\begin{enumerate}
\item[(i)]
$\dim \Hom(F_2/F_1,E_0)=l_2'$,
\item[(ii)]
$\dim \Ext^1(F_2/F_1,E_0)=0$ and
\item[(iii)]
$\psi$ is isomorphic in codimension 1.
\end{enumerate}
By the definition of Harder-Narasimhan filtration,
$\psi$ is not isomorphic. Thus $F_2/F_1$ is not locally free.
By (ii), $E^{\vee \vee} \cong F_1 \oplus F_2/F_1$,
which contradicts the locally freeness of $E^{\vee \vee}$.
Thus this case does not occur.
By (I-a), (I-b-1) and (I-b-2), we get a desired bound for the case (I).

Case (II).
We divide our consideration into three cases
 (II-a) $\langle v^2 \rangle>0$, 
 (II-b) $\langle v^2 \rangle=0$ and
 (II-c) $\langle v^2 \rangle<0$.

(II-a) If $\langle v_1^2 \rangle>0$, then
\begin{equation}
 \begin{split}
  \dim {\cal M}(v,n;v_1) &=\dim {\cal M}(v_1)+n(\rk v+1)\\
  &=\langle v_1^2 \rangle+1+n(\rk v+1)\\
  &=\langle v^2 \rangle+1-n(\rk v-1).
 \end{split}
\end{equation}
Hence if $\rk v \geq 3$, then $\codim_{{\cal M}(v)^s}
{\cal M}(v,n;v_1) \geq 2$. 
If $\rk v=2$, then the condition \eqref{eq:small} implies that
$a=3$. Hence we get $\langle v_1^2 \rangle=
\langle v^2 \rangle-2n \rk v \leq 0$.
Therefore this case does not occur.

(II-b) If $\langle v_1^2 \rangle=0$, then 
the argument in (I-b-2) implies that ${\cal M}(v_1)^{ss}$ consists of
non-locally free sheaves, which is a contradiction.

(II-c) We assume that $\langle v_1^2 \rangle<0$, that is,
$E^{\vee \vee}=E_0^{\oplus l}$. Then \eqref{eq:small} implies that 
$n \rk v_0-l \geq 0$.

(II-c-1) 
We first assume that $n \rk v_0-l \geq 1$. Then we get
\begin{equation}\label{eq:free}
 \begin{split}
  \codim_{{\cal M}(v)^s} {\cal M}(v,n) &=
  (2nl \rk v_0-2l^2+1)-(n(l\rk v_0+1)-l^2)\\
  &=n(l \rk v_0-1)-(l^2-1)\\
  & \geq \frac{l+1}{\rk v_0}(l \rk v_0-1)-(l^2-1)\\
  &=\frac{(l+1)(\rk v_0-1)}{\rk v_0} \geq 0,
 \end{split}
\end{equation}
and the equality holds if and only if
$\rk v_0=1$ and $n \rk v_0-l=1$.
By the computation of \eqref{eq:free}, it is easy to show that
$\codim_{{\cal M}(v)^s} {\cal M}(v,n) \geq 2$, if $\rk v_0 \geq 2$,
or $n \rk v_0-l \geq 2$.

(II-c-2) 
If $n \rk v_0-l =0$, then the primitivity of $v$ implies that
$n=1$.
Therefore we get a desired bound for the case (II).

By (I) and (II),
we complete the proof of Lemma \ref{lem:not-free}.
\qed

\subsubsection{Components containing locally free sheaves}

We shall prove Theorem \ref{thm:main} unless
$v=(\rk v_0)v_0-\omega$, or $\rk v_0=1$ and $v=l v_0-(l+1)\omega$. 
For a locally free sheaf $E \in {\cal M}(v)^s$,
we consider the dual of $E$.
Since $\langle v_0,v \rangle =a\rk v_0-2l<0$ and $E$ is stable,
$l':=\dim \Ext^2(E_0,E) \geq 2l-a\rk v_0>0$.
By Serre duality, we get an exact sequence
\begin{equation}\label{eq:reflection}
0 \to (E_0^{\vee})^{\oplus l'} \to E^{\vee} \to F \to 0,
\end{equation}
where $F \in {\cal M}(v^{\vee}-l'v_0^{\vee})^{\mu ss}$.
Then we see that
\begin{equation}\label{eq:9/21}
 \Hom(E_0^{\vee},F)=\Ext^2(E_0^{\vee},F)=0.
\end{equation}
Since $l'<l$, we see that
\begin{equation}
\begin{split}
\langle v(F)^2 \rangle &=
2l(a \rk v_0-l)-2l'(l'-2l+a \rk v_0)\\
&=2(l-l')(a \rk v_0-l+l')\\
& \geq 2(l-l')l>2(l-l')^2.
\end{split}
\end{equation}
Hence, by Lemma \ref{lem:quot}, we get
$\dim {\cal M}(v^{\vee}-l'v_0^{\vee})^{\mu ss}=\langle (v-l'v_0)^2 \rangle+1$.
Taking into account \eqref{eq:9/21}, we see that
the moduli number of $E^{\vee}$
which fits in the exact sequence \eqref {eq:reflection}
is given by
\begin{equation}\label{eq:brill}
 \dim{\cal M}(v^{\vee}-l'v_0^{\vee})^{\mu ss}+
 \dim Gr(\langle v^{\vee}-l'v_0^{\vee},v_0^{\vee} \rangle,l')
 =\langle v^2 \rangle+1-l'(l'+\langle v_0,v \rangle),
\end{equation}
(cf. [Y5, Lem. 2.6]).
We set
\begin{equation}
M_H(v)^*:=\left\{E \in M_H(v)^{loc}\left|
\begin{aligned}
& \Ext^1(E_0,E)=0,\\
&\coker(\Hom(E,E_0) \otimes E_0^{\vee} \to E^{\vee}) 
\in M_H(R_{v_0}(v^{\vee}))
\end{aligned}
\right. \right\}.
\end{equation}
Then, by using Lemma \ref{lem:quot}, Lemma \ref{lem:not-free} 
and \eqref{eq:brill}, we see that
$\codim_{M_H(v)}(M_H(v) \setminus M_H(v)^*) \geq 2$.
Since Theorem \ref{thm:main} holds for $M_H(R_{v_0}(v^{\vee}))$,
$M_H(R_{v_0}(v^{\vee}))$ is irreducible.
Therefore the morphism $M_H(v)^* \to M_H(R_{v_0}(v^{\vee}))$ is birational, 
if $M_H(v) \ne \emptyset$.
Conversely, for a $\mu$-stable vector bundle $F \in M_H(R_{v_0}(v^{\vee}))$,
we consider the universal extension
\begin{equation}\label{eq:univE'}
 0 \to \Ext^1 (F,E_0^{\vee})^{\vee} \otimes E_0^{\vee} \to E' \to F \to 0.
\end{equation}
We claim that $(E')^{\vee}$ is a stable sheaf
of $v((E')^{\vee})=v$.
This means that $M_H(v) \ne \emptyset$ and
Theorem \ref{thm:main} (1), (2-1) hold for this case.
The proof of Theorem \ref{thm:main} (2-2) is similar to
the proof for case A.
\newline
Proof of the claim:
Assume that $(E')^{\vee}$ is not stable.
Since $(E')^{\vee}$ is a $\mu$-semi-stable vector bundle,
there is a subbundle $G$ such that
\begin{enumerate}
\item[(1)]
$(E')^{\vee}/G$ is torsion free,
\item[(2)]
$v(G)=l_1 v_0- a_1 \omega$ and
\item[(3)]
$a_1/l_1<a/l$.
\end{enumerate}
Since $F^{\vee}$ is $\mu$-stable and $(E')^{\vee}/G$ is torsion free,
$F^{\vee} \cap G=0$, or $F^{\vee}$.
If $F^{\vee} \cap G=0$, then $G \to \Ext^1(F,E_0) \otimes E_0$ is injective.
Since the slopes of $G$ and $E_0$ are the same
and $G$ is locally free, we see that $G \cong E_0^{\oplus l_1}$,
which means that \eqref{eq:univE'} is not the universal extension.
Therefore $F^{\vee} \cap G$ must be equal to $F^{\vee}$.
This means that $G$ contains $F^{\vee}$ and 
$v(G/F^{\vee})=(l_1+l-a \rk v_0)v_0-(a_1-a) \omega$.
Since $G/F^{\vee}$ is a subsheaf of $\Ext^1(F,E_0) \otimes E_0$,
$a_1-a \geq 0$.
On the other hand, (3) and $l_1 \leq l$ implies that
$a_1<a$, which is a contradiction.
Hence $(E')^{\vee}$ is stable. 
\qed

\subsubsection{Non-locally free components}
We shall prove Theorem \ref{thm:main} for
\begin{equation}
 v=
  \begin{cases}
   (\rk v_0)v_0-\omega,\\
   lv_0-(l+1) \omega,\;\rk v_0=1.
  \end{cases}
\end{equation}
\begin{prop}\label{prop:FM}
Assume that $v=(\rk v_0)v_0-\omega$.
Then
$M_H(v) \cong X$.
\end{prop}
\begin{pf}
By the argument in (I-b-2),
$E \in M_H(v)$ satisfies that $E^{\vee \vee}=E_0^{\oplus \rk v_0}$.
Since $v=v(E_0^{\oplus \rk v_0})-\omega$, $E$ is the kernel of a quotient
$E_0^{\oplus \rk v_0} \to {\Bbb C}_x, x \in X$.
We shall construct a family of stable sheaves
$\{{\cal E}_x \}_{x \in X}$ of $v({\cal E}_x)=v$
and prove our proposition. 
For convenience sake, we set $X_i:=X$, $i=1,2$.
Let $\Delta \subset X_2 \times X_1$ be the diagonal.
We denote the projections $X_2 \times X_1 \to X_i$, $i=1,2$ by $p_i$.
We shall consider the evaluation map
\begin{equation}
 \phi:E_0^{\vee} \boxtimes E_0 \to (E_0^{\vee} \boxtimes E_0)_{|\Delta}
 \to {\cal O}_{\Delta}.
\end{equation}  
Then it is easy to see that
$\phi_x:=\phi_{|\{x \} \times X_1}, x \in X_2$ is surjective and
the induced homomorphism
\begin{equation}
 \Hom(E_0,(E_0^{\vee} \boxtimes E_0)_{|\{x \} \times X})
 \to \Hom(E_0,{\Bbb C}_x)
\end{equation}
is an isomorphism.
Hence $\ker \phi_x$ is stable.
Since ${\cal O}_{\Delta}$ is flat over $X_2$,
${\cal E}:=\ker \phi$ is flat over $X_2$ and
${\cal E}_{|\{x \} \times X}=\ker \phi_x, x \in X_2$ is stable.
Thus we get a morphism $X_2 \to M_H(v)$,
which is an isomorphism.
\end{pf}
\begin{cor}\label{cor:FMT}
$R_{v_0}$ is the Fourier-Mukai transform defined by ${\cal E}$.
\end{cor}
\begin{pf}
Let ${\cal F}:{\mathbf D}(X_1) \to {\mathbf D}(X_2)$ be the functor
defined by
\begin{equation}
{\cal F}(x):={\mathbf R}\Hom_{p_2}({\cal E},p_1^*(x)), x \in {\mathbf D}(X_1)
\end{equation}
and
$\widehat{{\cal F}}:{\mathbf D}(X_2) \to {\mathbf D}(X_1)$
the functor 
defined by 
\begin{equation}
\widehat{{\cal F}}(y):={\mathbf R}p_{1*}({\cal E} \otimes p_2^*(y)),
y \in {\mathbf D}(X_2).
\end{equation}
Then $\widehat{{\cal F}}[2]$ gives the inverse of ${\cal F}$.
Let $E$ be a coherent sheaf on $X_1$
such that
\begin{equation}\label{eq:cond}
 \begin{split}
  & \Hom(E_0,E)=0,\\
  & \Ext^2(E_0,E)=0.
 \end{split}
\end{equation}
We shall prove that $E$ satisfies ${\mathrm WIT}_1$ for ${\cal F}$, i.e.
\begin{equation}\label{eq:wit}
 \Hom_{p_2}({\cal E},p_1^*(E))=\Ext^2_{p_2}({\cal E},p_1^*(E))=0
\end{equation}
and ${\cal F}^1(E):=\Ext^1_{p_2}({\cal E},p_1^*(E))$ fits in 
the universal extension of $E$ by $E_0$:
Since ${\cal O}_{\Delta}$ is flat over $X_1$ and $K_X$ is trivial, we get
\begin{equation}\label{eq:local}
{\cal H}om^i({\cal O}_{\Delta},p_1^*(E))=
\begin{cases}
0, \;i=0,1,\\
{\cal O}_{\Delta} \otimes p_1^*(E), \;i=2.
\end{cases}
\end{equation}
By using local-global spectral sequence,
\eqref{eq:cond} and \eqref{eq:local}, we see that
\eqref{eq:wit} holds and  
we get an exact sequence
\begin{equation}\label{eq:univ}
 0 \to \Ext^1(E_0,E) \otimes E_0 \to
 \Ext^1_{p_2}({\cal E},p_1^*(E)) \to 
 E \to 0.
\end{equation}
We claim that this sequence gives the universal extension of
$E$ by $E_0$.
\newline
Proof of the claim:
If it is not the universal extension,
since $\Ext^1(E_0,E) \cong \Ext^1(E,E_0)^{\vee}$,
$E_0$ must be a direct summand of 
${\cal F}^1(E)$.
Since $E$ satisfies ${\mathrm WIT}_1$ for ${\cal F}$,
${\cal F}^1(E)$ satisfies ${\mathrm WIT}_1$ for $\widehat{{\cal F}}$,
which implies that $E_0$ also satisfies ${\mathrm WIT}_1$ for 
$\widehat{{\cal F}}$.
In particular, $R^2p_{1*}({\cal E} \otimes p_2^*(E_0))=0$.
On the other hand, a direct computation shows that
 $R^2p_{1*}({\cal E} \otimes p_2^*(E_0)) \cong E_0$,
which is a contradiction.
Therefore \eqref{eq:univ} is the universal extension. 
Thus ${\cal F}^1(E)$ is the reflection of $E$ by $v_0$.
\end{pf}

\begin{rem}
As in [Y7], we define
${\cal F}^H:H^*(X_1,{\Bbb Z}) \to H^*(X_2,{\Bbb Z})$ by
\begin{equation}
 {\cal F}^H(x):=p_{2*}\left(\ch({\cal E})^{\vee}
 p_1^*\sqrt{\td_{X_1}}p_2^*\sqrt{\td_{X_2}}
 p_{1}^*(x) \right), x \in H^*(X_1,{\Bbb Z}).
\end{equation}
Since $\ch({\cal E})^{\vee}=p_{1}^*(\ch(E_0))^{\vee}
p_{2}^*(\ch(E_0))-\ch({\cal O}_{\Delta})^{\vee}$, we also see that
${\cal F}^H(x)=-(x+\langle x,v(E_0) \rangle v(E_0))=-R_{v(E_0)}(x)$.
Thus we get the following commutative diagram:
\begin{equation}
 \begin{CD}
  {\mathbf D}(X_1) @>{\cal F}>> {\mathbf D}(X_2)\\
  @V{\ch \sqrt{\td_{X_1}}}VV @VV{\ch \sqrt{\td_{X_2}}}V\\
  H^*(X_1,{\Bbb Z}) @>>{-R_{v(E_0)}}> H^*(X_2,{\Bbb Z})
 \end{CD}
\end{equation}
\end{rem}

Finally we shall treat $M_H(l v_0-(l+1)\omega)$, $\rk v_0=1$.
\begin{prop}\label{prop:nlf2}
If $\rk v_0=1$, then 
$M_H(l v_0-(l+1)\omega) \cong \Hilb_X^{l+1}$ and $\theta_v$ is an isometry
of Hodge structures.
\end{prop}
\begin{pf}
We may assume that $v_0=1+\omega$.
Let $E$ be an element of $M_H(l-\omega)$.
We shall first prove that $E^{\vee \vee} \cong {\cal O}_X^{\oplus l}$.
To see this, it is sufficient to prove that
\begin{equation}\label{eq:pinch}
 n:=\dim(E^{\vee \vee}/E)=l+1.
\end{equation} 
Proof of \eqref{eq:pinch}:
Since $\chi(E)=l-1$ and $E$ is stable,
Serre duality implies that $\dim H^0(X,E^{\vee})\geq l-1$.
Hence we get an exact sequence
\begin{equation}\label{eq:pinch2}
 0 \to {\cal O}_X^{\oplus (l-1)} \to E^{\vee} \to I_Z \to 0
\end{equation}
where $I_Z \in \Hilb_X^{l+1-n}$.
If $n=0$, then $E$ is locally free.
By taking the dual of \eqref{eq:pinch2}, we get a
section of $E$, which contradicts the stability of $E$.
We assume that $0<n <l+1$.
Then $\dim \Ext^1(I_Z,{\cal O}_X)=\dim H^1(X,I_Z) =l-n$.
Hence we get a decomposition
$E^{\vee}={\cal O}_X^{\oplus (n-1)} \oplus F$,
where $F$ fits in an exact sequence
\begin{equation}
 0 \to {\cal O}_X^{\oplus (l-n)} \to F \to I_Z \to 0.
\end{equation}
Then $E^{\vee \vee}$ has a subsheaf ${\cal O}_X^{\oplus n}$.
The stability of $E$ implies that
\begin{equation}
 \frac{\chi({\cal O}_X^{\oplus n} \cap E)}n<
 \frac{\chi(E)}{l}.
\end{equation}
Since $\chi({\cal O}_X^{\oplus n} \cap E) \geq 
\chi({\cal O}_X^{\oplus n})-n=n$, this is impossible.
Therefore \eqref{eq:pinch2} holds, 
which implies that $E^{\vee \vee}={\cal O}_X^{\oplus l}$.

Conversely for a general quotient 
$\phi:{\cal O}_X^{\oplus l} \to \oplus_{i=1}^{l+1}{\Bbb C}_{x_i},
x_1,x_2,\dots,x_{l+1} \in X$,
it is easy to see that $\ker \phi$ is stable. 
Thus $M_H(l-\omega) \ne \emptyset$.

We shall prove that it is isomorphic to $\Hilb_X^{l+1}$.
For this purpose, we shall consider a functor 
${\cal G}:{\bf D}(X_1) \to {\bf D}(X_2)_{op}$ which is the composition of 
reflection by $v({\cal O}_X)$ with the taking dual functor:
\begin{equation}
 {\cal G}(x):={\bf R}\Hom_{p_2}(p_1^*(x),I_{\Delta}),
 x \in {\bf D}(X_1),
\end{equation}
where we use the same notation as in Proposition \ref{prop:FM} 
and ${\bf D}(X_2)_{op}$ is the opposite category of ${\bf D}(X_2)$.
Then ${\cal G}$ gives an equivalence of categories.
For $E \in M_H(l-\omega)$, we shall prove that 
\begin{enumerate}
\item[(a)] 
 $\Ext^i_{p_2}(p_1^*(E),I_{\Delta})=0$, $i=0,2$ and
\item[(b)]
 ${\cal G}^1(E):=\Ext^1_{p_2}(p_1^*(E),I_{\Delta})$ is an ideal sheaf of
 colength $l+1$.
\end{enumerate}
Then the map $M_H(l-\omega) \to \Hilb_X^{l+1}$ sending
$E$ to ${\cal G}^1(E)$ gives an isomorphism of moduli spaces.
The second assertion follow from [Y7, Prop. 2.5] or a direct computation
by using the equality
\begin{equation}
 {\cal G}({\cal E})={\mathbf R}\Hom_{p_{M_H(l-\omega)}}
 ({\cal E},{\cal O}_{M_H(l-\omega) \times X_1})\boxtimes {\cal O}_{X_2}
 -{\mathbf R}{\cal H}om({\cal E},{\cal O}_{M_H(l-\omega) \times X_2})
\end{equation}
as an element of Grothendieck group of $M_H(l-\omega) \times X_2$,
where ${\cal E}$ is a quasi-universal family on $M_H(l-\omega) \times X$.

Proof of (a), (b):
By Serre duality and the stability of $E$,
$\Ext^2(E,I_x)=0$ for all $x \in X$.
Also we see that $\Hom(E,I_x)=0$ for $x \not \in 
\Supp(E^{\vee \vee}/E)$.
Hence by the base change theorem and its proof, we see that
(a) holds and ${\cal G}^1(E)$ is torsion free.
It is easy to see that $v({\cal G}(E))=-R_{v({\cal O}_X)}(v(E)^{\vee})
=-1+l\omega$.
Hence $v({\cal G}^1(E))=1-l \omega$. Therefore 
${\cal G}^1(E)$ is an ideal sheaf of colength $l+1$.

\end{pf}

\subsection{The case of $\langle v^2 \rangle=2l^2$}

We next treat the case where $\langle v^2 \rangle=2l^2$.
Let $E$ be a general member of $M_H(v) \setminus M_H(v)^{\mu s}$.
Then by the proof of Lemma \ref{lem:quot},
$E$ fits in an exact sequence
\begin{equation}
0 \to E_1 \to E \to E_0 \to 0
\end{equation}
where $E_1$ is a $\mu$-stable vector bundle.

Indeed, if $\langle v^2 \rangle=2l^2$, then
the primitivity of $v$ implies that (i) $v=lv_0-2\omega$, $l=\rk v_0$
or (ii) $v=lv_0-\omega$, $2l=\rk v_0$. In particular $\rk v_0 >1$. 
In the notation of the proof of Lemma \ref{lem:quot}, we see that
$k=1$. Hence $E$ fits in the above exact sequence.
Then in the same way as in 2.2.2,
we get Theorem \ref{thm:main} in this case.   
\qed

\subsection{Some remarks}
By the proof of Theorem \ref{thm:main} for case B and Lemma \ref{lem:quot},
we also get the following.
\begin{enumerate}
\item[(1)]
$M_H(v)^{\mu s} \ne \emptyset$ if and only if $\langle v^2 \rangle \geq 2l^2$.
\item[(2)]
$M_H(v)^{loc}=\emptyset$ if and only if (i) $\rk v=1$, 
(ii) $v=(\rk v_0)v_0-\omega$, or
(iii) $\rk v_0=1$ and $v=lv_0-(l+1) \omega$.
\end{enumerate}
Combining Propositions \ref{prop:nlf2}, \ref{prop:FM}
and Remark \ref{rem:A},
we get Proposition \ref{prop:nlf}.

\section{Relation to Montonen-Olive duality}


In this section, we shall consider the relation between Theorem
\ref{thm:main} and Montonen-Olive duality in Physics.
Roughly speaking, Montonen-Olive duality says that
the generating function of Euler characteristics of moduli spaces of
vector bundles becomes a modular form.
In this paper, 
we concentrate on moduli spaces of vector bundles on K3 surfaces.
We first describe physical predictions and their modifications.
For more details and related results, see [MNVW], [V-W] and 
[G\"{o}3], [Y4], [Y6].

\subsection{Physical predictions}
We fix a K3 surface $X$.
We regard $H^2(X,{\Bbb Z})$ as a lattice by a bilinear form
$Q(x,y)=-\int_X xy, x,y \in H^2(X,{\Bbb Z})$.
Let $P$ be a orthogonal decomposition of
$H^2(X,{\Bbb Z}) \otimes {\Bbb R}$ 
as a sum of definite signature:
\begin{equation}
P:H^2(X,{\Bbb Z}) \otimes {\Bbb R} \cong {\Bbb R}^{19,0} \oplus
 {\Bbb R}^{0,3}.
\end{equation}
Let $P_L(x)=x_L$, $P_R(x)=x_R$ denote the projections onto the two factors.

For $v=r+\xi+a\omega \in H^*(X,{\Bbb Z})$ of $\xi \in H^2(X,{\Bbb Z})$,
we choose a suitable complex structure such that
$\xi$ become holomorphic.
Then we define $M(v)$ as a moduli space of stable sheaves 
on this surface.
Let $Z_r(\tau,x)$ be $U(r)$-partition function defined in
[MNVW, sect. 3]:
\begin{equation}\label{eq:u(r)}
 Z_r(\tau,x):=\sum
 \begin{Sb}
 v \in H^*(X,{\Bbb Z})\\
 \rk v=r
 \end{Sb}
 ``\chi(M(v))"q^{\frac{\langle v^2 \rangle}{2r}} 
 q^{\frac{1}{2r}Q(c_1(v)_L^2)}
 \overline{q}^{\frac{-1}{2r}Q(c_1(v)_R^2)} e^{Q(c_1(v),x)}
\end{equation} 
where 
$(\tau, x) \in {\Bbb H} \times H^2(X,{\Bbb Z}) \otimes {\Bbb C}$,
$q:=\exp(2 \pi \sqrt{-1} \tau)$, $e:=\exp(2 \pi \sqrt{-1})$ and  
$``\chi(M(v))"$ is a kind of ``Euler characteristics" of
a nice compactification of $M(v)$.

\begin{rem}
More precisely, Minahan et al. considered $Z_r(\tau,0)$.
Combining the computations in [MNVW, sect. 6], 
we propose the definition \eqref{eq:u(r)}.
\end{rem}
Unfortunately, there is no mathematical definition of 
this ``Euler characteristics".
Since $M(v)$ is smooth and compact for primitive $v$,
we can expect that $``\chi(M(v))"$ coincides with the 
ordinary Euler characteristics $\chi(M(v))$.

Then Montonen-Olive duality for $U(r)$ gauge group
asserts that 
\begin{itemize}
\item[($\#$)]
$Z_r(\tau,x)$ transforms like a Jacobi form
of holomorphic/anti-holomorphic weight 
$$(-\chi(X)/2+b_{-}(X)/2,b_{+}(X)/2)=(-5/2,3/2)$$ 
\end{itemize}
(cf. [E-Z]).
For $\alpha \in H^2(X,{\Bbb Z})$,
let $Z_r^{\alpha}(\tau)$ be $PSU(r)$-partition function defined in
[V-W]:
\begin{equation}
 Z_r^{\alpha}(\tau):=\sum
 \begin{Sb}
 v \in H^*(X,{\Bbb Z})\\
 \rk v=r, c_1(v)=\alpha 
 \end{Sb}
 ``\chi(M(v))"q^{\frac{\langle v^2 \rangle}{2r}}. 
\end{equation}
Then 
\begin{equation}\label{eq:factor}
Z_r(\tau,x)=\sum_{\alpha \in H^2(X,{\Bbb Z})/rH^2(X,{\Bbb Z})}
Z_r^{\alpha}(\tau) \Theta_{\alpha,r}(\tau,P,x),
\end{equation}
where 
\begin{equation}
\Theta_{\alpha,r}(\tau,P,x)=\sum_{c \in \alpha+rH^2(X,{\Bbb Z})}
q^{\frac{1}{2r}Q(c_L^2)}
 \overline{q}^{\frac{-1}{2r}Q(c_R^2)} e^{Q(c,x)}
\end{equation}
is Siegel-Narain theta function (cf. [M-W, Appendix B]).
If $r=1$, then it is known that 
$Z_1^0(\tau)=\frac{1}{\eta(\tau)^{24}}$ ([G\"{o}1], [V-W]).
Hence 
\begin{equation}
\begin{split}
Z_1(\tau,x) &=Z_1^0(\tau)
\left(\sum_{c \in H^2(X,{\Bbb Z})} q^{\frac{1}{2}Q(c_L^2)}
 \overline{q}^{\frac{-1}{2}Q(c_R^2)} e^{Q(c,x)} \right)\\
 &= \frac{1}{\eta(\tau)^{24}}\Theta(\tau,P,x),
\end{split}
\end{equation}
where $\Theta(\tau,P,x)= \sum_{c \in H^2(X,{\Bbb Z})} 
q^{\frac{1}{2}Q(c_L^2)}
 \overline{q}^{\frac{-1}{2}Q(c_R^2)} e^{Q(c,x)} $.
Since $\Theta(\tau,P,x)$ transforms like a Jacobi form of 
holomorphic/anti-holomorphic
weight $(19/2,3/2)$,
$Z_1(\tau,x)$ transforms like a Jacobi form of holomorphic/anti-holomorphic
weight $(-5/2,3/2)$:
\begin{equation}
Z_1\left(-\frac{1}{\tau},\frac{x_L}{\tau}+\frac{x_R}{\overline{\tau}}\right)
=(-\sqrt{-1}\tau)^{-5/2}(\sqrt{-1}\overline{\tau})^{3/2}
e^{\frac{Q(x_L^2)}{2\tau}}
e^{\frac{Q(x_R^2)}{2\overline{\tau}}}Z_1(\tau,x_L+x_R).
\end{equation}
Then $Z_r(\tau,x)$ is given by Hecke transformation of 
order $r$ of $Z_1(\tau,x)$ ([MNVW]):
\begin{equation}\label{eq:U(r)}
Z_r(\tau,x)=\frac{1}{r^2} 
\sum
\begin{Sb}
a,b,d \geq 0\\
ad=r\\
b<d
\end{Sb}
d Z_1\left(\frac{a \tau+b}{d},ax \right).
\end{equation}
In particular, $Z_r(\tau,x)$ transforms like a Jacobi form of
holomorphic/anti-holomorphic
weight $(-5/2,3/2)$ and index $r$.
Thus ($\#$) holds.
\begin{rem}
For $PSU(r)$-partition functions, we get
the following:
\begin{equation}
Z_r^{c_1}(\tau)=
\frac{1}{r^2} 
\sum
\begin{Sb}
a,b,d \geq 0\\
ad=r\\
b<d\\
a \xi=c_1
\end{Sb}
d Z_1^0\left(\frac{a \tau+b}{d}\right)
e^{(-\frac{b}{2d}(\xi^2))}.
\end{equation}
Combining \eqref{eq:factor} with the transformation law
\begin{equation}
\begin{split}
\Theta_{\alpha,r}\left(\frac{-1}{\tau},P,\frac{x_L}{\tau}+
\frac{x_R}{\overline{\tau}}\right)
&=r^{-11}(-\sqrt{-1}\tau)^{19/2}(\sqrt{-1}\overline{\tau})^{3/2}
e^{\frac{rQ(x_L^2)}{2\tau}}
e^{\frac{rQ(x_R^2)}{2\overline{\tau}}}\\
 & \quad \quad \cdot\left(\sum_{\beta \in H^2(X,{\Bbb Z}/r{\Bbb Z})}
e^{\frac{-Q(\alpha,\beta)}{r}}\Theta_{\beta,r}(\tau,x_L+x_R)\right),
\end{split}
\end{equation}
we can deduce from ($\#$) the following transformation law:
\begin{equation}
Z_r^{\alpha}(-1/\tau)=
r^{-11}(-\sqrt{-1}\tau)^{-12}
\sum_{\beta \in H^2(X,{\Bbb Z}/r{\Bbb Z})}
e^{\frac{Q(\alpha,\beta)}{r}}Z_r^{\beta}(\tau).
\end{equation}
This formula is of course compatible with 
Montonen-Olive duality for $PSU(r)$ group [V-W].
\end{rem}

\subsection{Relation to Theorem \ref{thm:main}}
We shall check that Theorem \ref{thm:main} is compatible with
\eqref{eq:U(r)}.
For simplicity, we set $X^{[n]}=\Hilb_X^n$.

\begin{equation}
\begin{split}
\sum_{0 \leq b<d}d Z_1\left(\frac{a \tau+b}{d},ax \right) &=
\sum_{0 \leq b<d}\sum_{\xi \in H^2(X,{\Bbb Z})}
\sum_n d\chi(X^{[n]})q^{\frac{a}{d}(n-1)}
 q^{\frac{a}{2d}Q(\xi_L^2)}
 \overline{q}^{\frac{-a}{2d}Q(\xi_R^2)} e^{aQ(\xi,x)} 
 e^{\frac{b}{d}((n-1)+Q(\xi^2)/2)}\\
 &=\sum_{\xi \in H^2(X,{\Bbb Z})}
\sum_{d|n-1+\frac{Q(\xi^2)}{2}}
 d^2\chi(X^{[n]})q^{\frac{a}{d}(n-1)}
 q^{\frac{a}{2d}Q(\xi_L^2)}
 \overline{q}^{\frac{-a}{2d}Q(\xi_R^2)} e^{aQ(\xi,x)} \\
 &=\sum_{\xi \in H^2(X,{\Bbb Z})}
\sum_{k}
 d^2\chi(X^{[kd-Q(\xi^2)/2+1]})q^{\frac{a}{d}(kd-Q(\xi^2)/2)}
 q^{\frac{a}{2d}Q(\xi_L^2)}
 \overline{q}^{\frac{-a}{2d}Q(\xi_R^2)} e^{aQ(\xi,x)} \\
&= \sum_{\xi \in H^2(X,{\Bbb Z})}
\sum_{w=(d,\xi,-k)}
 d^2\chi(X^{[\langle w^2 \rangle/2+1]})
 q^{\frac{a}{2d}\langle w^2 \rangle}
 q^{\frac{a}{2d}Q(c_1(w)_L^2)}
 \overline{q}^{\frac{-a}{2d}Q(c_1(w)_R^2)} e^{aQ(c_1(w),x)} \\
&= \sum_{\rk w=d}
 d^2\chi(X^{[\langle w^2 \rangle/2+1]})
 q^{\frac{1}{2r}\langle (aw)^2 \rangle}
 q^{\frac{1}{2r}Q(c_1(aw)_L^2)}
 \overline{q}^{\frac{-1}{2r}Q(c_1(aw)_R^2)} e^{Q(c_1(aw),x)}.
 \end{split}
 \end{equation}
Therefore we get 
\begin{equation}
 ``\chi(M(v))"=\sum_{v=aw}\frac{1}{a^2}
 \chi(X^{[\langle w^2 \rangle/2+1]}).
 \end{equation}
 In particular, if $v$ is primitive, then by Corollary \ref{cor:hodge},
we get
 \begin{equation}
 ``\chi(M(v))"=\chi(X^{[\langle v^2 \rangle/2+1]})=\chi(M(v)).
 \end{equation}
 This implies that $\chi(M(v))$ is related to modular forms  
 and in particular Hecke transforms. 
To understand the meaning of $``\chi(M(v))"$ for non-primitive $v$
is a challenging problem.
The relation to O'Grady's symplectic compactification of
$M(2-2\omega)$ ([O2]) is also an interesting problem.

\section{Appendix}

In this appendix, we shall explain our method for dimension counting
of substacks of ${\cal M}(v)^{\mu ss}$.
Since most results are appeared in another forms
(cf. [D-R], [H-N]),
we only give an outline.

\subsection{Notation.}
For an ample divisor $H'$ on $X$,
let $Q(mH',v)$ be the open subscheme of 
the quot-scheme $\Quot_{{\cal O}_X(-mH')^{\oplus N}/X/{\Bbb C}}$
consisting of points
\begin{equation}
\lambda:{\cal O}_X(-mH')^{\oplus N} \to E
\end{equation}
such that
\begin{enumerate}
\item
$v(E)=v$,
\item
$\lambda$ induces an isomorphism
$H^0(X,{\cal O}_X^{\oplus N}) \cong H^0(X,E(mH'))$,
\item
$H^i(X,E(mH'))=0$, $i>0$.
\end{enumerate}
Let ${\cal O}_{Q(mH',v) \times X}(-mH')^{\oplus N} \to {\cal Q}_v$
be the universal quotient.
We set $V_v:={\cal O}_X(-mH')^{\oplus N}$.
For our purpose, the choice of $mH$ is not so important.
Hence we simply denote $Q(mH',v)$ by $Q(v)$.

Let $q_v:Q(v) \to {\cal M}(v)$ be the natural map.
We denote the pull-backs $q_v^{-1}({\cal M}(v)^{\mu ss}),
 q_v^{-1}({\cal M}(v)^{ss}),\ldots$ by
$Q(v)^{\mu ss}, Q(v)^{ss}, \dots $ respectively.
If we choose a suitable $Q(v)$, then
$q_v:Q(v)^{\mu ss} \to {\cal M}(v)^{\mu ss}$ is 
surjective and 
${\cal M}(v)^{\mu ss}$ is a quotient stack of $Q(v)^{\mu ss}$
by a natural action of $G_v:=GL(N)$:
\begin{equation}
 {\cal M}(v)^{\mu ss} \cong \left[Q(v)^{\mu ss}/G_v \right].
\end{equation}
From now on, we assume that $q_v:Q(v)^{\mu ss} \to {\cal M}(v)^{\mu ss}$
is surjective.

\subsection{Stack of filtrations}

\begin{defn}
${\cal F}(v_1,v_2)$ is the stack of filtrations
$F_1 \subset E$, $E \in {\cal M}(v)$ such that
\begin{enumerate}
\item
$F_1$ is a $\mu$-semi-stable sheaf of $v(F_1)=v_1$.
\item
$E/F_1$ is a $\mu$-semi-stable sheaf of $v(E/F_1)=v_2$. 
\end{enumerate}
Let $p_v:{\cal F}(v_1,v_2) \to {\cal M}(v)^{\mu ss}$
be the projection sending $(F_1 \subset E)$ to
$E$ and
$p_{v_1,v_2}:{\cal F}(v_1,v_2) \to {\cal M}(v_1)^{\mu ss} \times
{\cal M}(v_2)^{\mu ss}$ the morphism sending 
$(F_1 \subset E)$ to $(F_1,E/F_1)$.
\end{defn}
We consider an open subscheme $F(v_1,v_2)$ of 
$\Quot_{{\cal Q}_v/Q(v)^{\mu ss} \times X/Q(v)^{\mu ss}}$
consisting of quotients $({\cal Q}_v)_x \to E_2, x \in Q(v)^{\mu ss}$
such that
$E_2$ is a $\mu$-semi-stable sheaf of $v(E_2)=v_2$.
Then 
\begin{equation}
{\cal F}(v_1,v_2)=[F(v_1,v_2)/G_v].
\end{equation}
We shall give another expression of ${\cal F}(v_1,v_2)$
which is useful to compute the dimensions of substacks of
${\cal F}(v_1,v_2)$ and its projections to ${\cal M}(v)$.  

We shall choose $Q(mH',v)^{\mu ss}$, $Q(mH',v_1)^{\mu ss}$
and $Q(mH',v_2)^{\mu ss}$ for the same $mH'$.
Then $V_v=V_{v_1} \oplus V_{v_2}$. 
For simplicity, we set ${\cal Q}_i:={\cal Q}_{v_i}, V_i:=V_{v_i}, \ldots$
and ${\cal K}_i$ are the universal subsheaves of 
${\cal O}_{Q(v_i)^{\mu ss}} \otimes V_i$, $i=1,2$. 
We define a scheme $\varpi:Y \to Q(v_1)^{\mu ss} \times Q(v_2)^{\mu ss}$
by
\begin{equation}
 Y:=\{\psi:({\cal K}_2)_{x_2} \to ({\cal Q}_1)_{x_1}|
 (x_1,x_2) \in Q(v_1)^{\mu ss} \times Q(v_2)^{\mu ss} \}.
\end{equation}
Then $Y$ parameterizes subsheaves $K \subset V$
such that $K \cap V_1=({\cal K}_1)_{x_1}$ and
$K/K \cap V_1=({\cal K}_2)_{x_2}$:
For a quotient $\psi:({\cal K}_2)_{x_2} \to ({\cal Q}_1)_{x_1}$,
the subsheaf $K$ of $V_v$ is defined by
\begin{equation}
 K:=\{(a_1,a_2) \in V_1 \oplus V_2|a_2 \in ({\cal K}_2)_{x_2},
 \psi(a_2)=a_1 {\mathrm mod} ({\cal K}_1)_{x_1} \}.
\end{equation}
Considering the quotient
$V_v \to V_v/K$,
$Y$ also parameterizes the following exact and commutative diagram:
\begin{equation}
 \begin{CD}
  @. 0 @.0 @.0 @.\\
  @.@AAA @AAA @AAA @.\\
  0 @>>> ({\cal Q}_1)_{x_1} @>>> E @>>> ({\cal Q}_2)_{x_2} @>>>0\\
  @.@AAA @AAA @AAA @.\\
  0 @>>> V_1 @>>> V_v @>>> V_2 @>>>0
 \end{CD}
\end{equation}
Let $\xi:Y \times G_v \to F(v_1,v_2)$ be the morphism
sending $(y:V_v \to E, g) \in Y \times G_2$ to
\begin{equation}
(y \circ g:V_v \to E, E \to {\cal Q}_{x_2}) \in F(v_1,v_2),
\end{equation}
where $\varpi(y)=(x_1, x_2)$.
Let $P$ be the parabolic subgroup of $G_v$ fixing $V_1$.
Then $Y$ has a natural action of $P$ and
$\xi$ induces a morphism
$Y \times_P G_v \to F(v_1,v_2)$, which is $G_v$-equivariant.
It is easy to see that this morphism is an isomorphism
(cf. [Y2, appendix]).
Therefore
\begin{equation}\label{eq:A-isom}
 \begin{split}
  {\cal F}(v_1,v_2)& \cong\left[Y \times_P G_v/G_v \right]\\
   & \cong \left[Y/P\right].
 \end{split}
\end{equation}
By using \eqref{eq:A-isom}, we shall prove the following.
\begin{lem}\label{lem:A-est1}
We set
\begin{equation}
 \begin{split}
  {\cal N}^n(v_1,v_2):&=\{(E_1,E_2) \in {\cal M}(v_1)^{\mu ss} \times 
  {\cal M}(v_2)^{\mu ss}|
  \dim \Hom(E_1,E_2)=n \},\\
  {\cal F}^n(v_1,v_2):&=p_{v_1,v_2}^{-1}({\cal N}^n(v_1,v_2))\\
  &=\{(F_1 \subset E) \in {\cal F}(v_1,v_2)| \dim \Hom(F_1,E/F_1)=n \}.
 \end{split}
\end{equation}
Then, 
\begin{equation}
 \dim {\cal F}^n(v_1,v_2)=\dim {\cal N}^n(v_1,v_2)+
 \langle v_1,v_2 \rangle+n.
\end{equation} 
\end{lem}

\begin{pf}
We set 
\begin{equation}
Q^n(v_1,v_2):=\{(x_1,x_2) \in Q(v_1)^{\mu ss} \times Q(v_2)^{\mu ss}
|\dim \Hom(({\cal Q}_1)_{x_1},({\cal Q}_2)_{x_2})
=n \}.
\end{equation}
For $(x_1, x_2) \in Q^n(v_1,v_2)$,
there is an exact sequence
\begin{equation}
 0 \to \Hom(({\cal Q}_2)_{x_2},({\cal Q}_1)_{x_1}) \to 
 \Hom(V_2,({\cal Q}_1)_{x_1}) \to 
 \Hom(({\cal K}_2)_{x_2},({\cal Q}_1)_{x_1}) \to
 \Ext^1(({\cal Q}_2)_{x_2},({\cal Q}_1)_{x_1}) \to 0.
\end{equation}
Since $\dim \Hom (V_2,({\cal Q}_1)_{x_1})=\rk V_1 \rk V_2$ and
$\dim \Ext^2(({\cal Q}_2)_{x_2},({\cal Q}_1)_{x_1})=n$,
$\Hom({\cal K}'_2,{\cal Q}'_1)$ is a locally free sheaf of
rank $\langle v_1,v_2 \rangle+ \rk V_1 \rk V_2+n$ on $Q^n(v_1,v_2)$,
where ${\cal K}'_2$ and ${\cal Q}'_1$ are pull-backs of 
${\cal K}_2$ and ${\cal Q}_1$ to $Q^n(v_1,v_2)$ respectively.
We set $Y^n:={\Bbb V}(\Hom({\cal K}'_2,{\cal Q}'_1)^{\vee}) \to Q^n(v_1,v_2)$.
Then ${\cal F}^n(v_1,v_2) \cong [Y^n/P]$.
Hence we get that
\begin{equation}
 \begin{split}
  \dim {\cal F}^n(v_1,v_2)&=\dim Y^n-\dim P\\
  &=\dim Q^n(v_1,v_2)+\langle v_1,v_2 \rangle+n-
  (\dim G_1+\dim G_2)\\
  &=\dim {\cal N}^n(v_1,v_2)+\langle v_1,v_2 \rangle+n.
 \end{split}
\end{equation}
\end{pf} 
 By similar method as in the proof of Lemma \ref{lem:A-est1},
we can also prove the following. 

\begin{lem}\label{lem:A-est2}
Let ${\cal F}^0(v_1,v_2,\dots,v_s)$ be the stack of filtrations
\begin{equation}
0 \subset F_1 \subset F_2 \subset \dots \subset F_s=E, E \in {\cal M}(v)
\end{equation}
such that 
\begin{enumerate}
\item
 $F_i/F_{i-1}$, $1 \leq i \leq s$ are semi-stable of $v(F_i/F_{i-1})=v_i$.
\item
 $\Hom(F_i/F_{i-1},F_j/F_{j-1})=0$, $i<j$.
\end{enumerate}
Then
\begin{equation}
 \dim {\cal F}^0(v_1,v_2,\dots,v_s)=\sum_{i=1}^s {\cal M}(v_i)^{ss}
 +\sum_{i<j} \langle v_i,v_j \rangle.
\end{equation}
\end{lem}

\subsection{Supplement for the proof of Lemma \ref{lem:ss}}


We shall explain how to derive \eqref{eq:4.3} from \eqref{eq:j}.
Let $q_l:Q(lw)^{ss} \to \overline{M}_H(lw):=Q(lw)^{ss}/G_{lw}$ be the
quotient map.
For a sequence of positive integers 
$l_1 \leq l_2 \leq \dots \leq l_s$ of $\sum_{i=1}^s l_i=l$,
we set
\begin{equation}
 \begin{split}
  \overline{M}_H(lw;l_1,l_2,\dots,l_s):&=
  \{\oplus_{i=1}^s E_i^{\oplus l_i} \in \overline{M}_H(lw )|
  E_1,E_2,\dots, E_s \in M_H(w), \text{$E_i \ne E_j$ for $i \ne j$}\},\\
  Q(l w;l_1,l_2,\dots,l_s)^{ss} :&=
  q_l^{-1}(\overline{M}_H(lw;l_1,l_2,\dots,l_s)).
 \end{split}
\end{equation}
For simplicity, we set $G:=G_{lw}$ and $G_i:=G_{l_iw}$, $i=1,2,\dots,s$.
For quotients $\phi_i:V_i \to {\cal Q}_{x_i} \in Q(l_i w; l_i)^{ss}$,
$i=1,2,\dots,s$ and an element $g \in G$,
we define a quotient
\begin{equation}
 (\oplus_{i=1}^s \phi_i) \circ g:V \to \oplus_{i=1}^s {\cal Q}_{x_i}.
\end{equation}
It will define a morphism
$\prod_{i=1}^s Q(l_i w; l_i)^{ss} \times G \to Q(l w)$.
Let $\prod_{i=1}^s Q(l_i w; l_i)^{ss} \times_{\prod_i G_i}G$
be the quotient of $\prod_{i=1}^s Q(l_i w; l_i)^{ss} \times G$
by a natural action of ${\prod_i G_i}$.
Then the above morphism induces a morphism
\begin{equation}
 \pi_{l_1,l_2,\dots,l_s}:
 \prod_{i=1}^s Q(l_i w; l_i)^{ss} \times_{\prod_i G_i}G \to Q(lw).
\end{equation}
By the construction of $\pi_{l_1,l_2,\dots,l_s}$,
$\im \pi_{l_1,l_2,\dots,l_s}$ contains $Q(l w;l_1,l_2,\dots,l_s)^{ss}$.
Hence
\begin{equation}
\dim \left[Q(l w;l_1,l_2,\dots,l_s)^{ss}/G \right]
\leq \sum_{i=1}^s \dim \left[Q(l_i w;l_i)^{ss}/G_i \right].
\end{equation}
We shall prove that 
\begin{equation}\label{eq:Q(l_i w;l_i)/G_i}
\dim \left[Q(l_i w;l_i)^{ss}/G_i \right] \leq 1.
\end{equation}
Then $\dim \left[Q(l w;l_1,l_2,\dots,l_s)^{ss}/G \right] \leq s \leq l$.
Clearly
$\dim \left[Q(l w;1,1,\dots,1)^{ss}/G \right]=l$.
Hence $\dim[Q(lw)^{ss}/G]=l$.

Proof of \eqref{eq:Q(l_i w;l_i)/G_i}:
Since $q_{l_i}:Q(l_i w;l_i)^{ss} \to \overline{M}_H(l_i w;l_i)$
is surjective and $\dim \overline{M}_H(l_i w;l_i)=2$,
it is sufficient to prove that 
\begin{equation}
\dim \left[q_{l_i}^{-1}(E^{\oplus l_i})/G_i \right] \leq -1,\; E \in M_H(w).
\end{equation}
By definition, 
${\cal J}(l_i,E)=[q_{l_i}^{-1}(E^{\oplus l_i})/G_i]$.
Hence we obtain this claim from \eqref{eq:j}.

\vspace{1pc}

{\it Acknowledgement.}
This paper heavily depends on Mukai's wonderful works
[Mu1,2,3]. 
I would like to thank Shigeru Mukai 
for valuable discussion on his works.
I would also like to thank Toshiya Kawai
for useful discussions on Jacobi forms.
Main part of this work 
was done when I stayed at Max Planck Institut f\"{u}r Mathematik.
I would like to thank Max Planck Institut f\"{u}r Mathematik
for support and hospitality.

\end{document}